\documentclass{amsart}   
\usepackage{amssymb}
\usepackage{amscd}
\newcommand\alge{\mathcal A}

\newcommand\maA{\mathcal A}
\newcommand\maF{\mathcal F}

\newcommand\maC{\mathcal C}

\newcommand\maP{\mathcal P}

\newcommand\Ker{\operatorname{Ker}}

\newcommand\EH{\operatorname{EH}}
\newcommand\codim{\operatorname{codim}}

\newcommand\CAA{{\mathcal A}(M,F)}

\newcommand\CC{\mathbb C}

\newcommand\NN{\mathbb N}
\newcommand\RR{\mathbb R}
\newcommand\ZZ{\mathbb Z}

\renewcommand\S{\mathbb S}
\newcommand\pa{\partial}
\newcommand\Alg{\mathcal{A}}

\newcommand\GR{\mathcal{G}}
\newcommand\CI{\mathcal{C}^\infty}

\newcommand\HH{\operatorname{HH}}
\newcommand\Hd{\operatorname{HH}}

\newcommand\Hp{\operatorname{HP}}

\newcommand{\R}{\RR}

\newcommand{\C}{\CC}
\newcommand{\Z}{\ZZ}

\newcommand{\cohom}{\operatorname{H}}
\newtheorem{theorem}{Theorem}
\newtheorem{proposition}{Proposition}
\newtheorem{corollary}{Corollary}
\newtheorem{lemma}{Lemma}

\theoremstyle{definition}
\newtheorem{definition}{Definition}

\theoremstyle{remark}
\newtheorem{remark}{Remark}
\newtheorem{example}[theorem]{Example} 

\begin{document}

\title[Residues for foliations]
{Residues and homology for pseudodifferential
operators on foliations}

\author[M-T. Benameur]{Moulay-Tahar Benameur}
\address{Inst. Desargues, Lyon, France}
\email{benameur@desargues.univ-lyon1.fr}
\author[V. Nistor]{Victor Nistor}
\address{Pennsylvania State University,
University Park, PA 16802}
\email{nistor@math.psu.edu}

\thanks{Nistor was partially supported by NSF Grant DMS-9971951 and 
``collaborative research'' grant 9981251. Manuscripts available from {\bf 
http:{\scriptsize//}www.math.psu.edu{\scriptsize/}nistor{\scriptsize/}}.} 


\begin{abstract} We study the Hochschild homology groups of the algebra of
complete symbols on a foliated manifold $(M,F)$. The first step is to
relate these groups to the Poisson homology of $(M,F)$ and of other
related foliated manifolds. We then establish several general
properties of the Poisson homology groups of foliated manifolds. As
an example, we completely determine these Hochschild homology groups
for the algebra of complete symbols on the irrational slope foliation
of a torus (under some diophantine approximation assumptions). We also
use our calculations to determine all residue traces on algebras of
pseudodifferential operators along the leaves of a foliation.
\end{abstract}
\maketitle
\tableofcontents

\section*{Introduction\label{Sect.I}}

This paper is a continuation of \cite{BenameurNistor1} and
\cite{BenameurNistor2}.  In those papers we have determined the
Hochschild, cyclic, and periodic cyclic homology of certain algebras of
complete symbols defined using groupoids. Our results were complete
for periodic cyclic homology, these groups being given directly in
terms of the cohomology of the cosphere bundle of the algebroid
associated to our groupoid (this result is recalled in Theorem
\ref{theorem.prev}), for any differentiable groupoid.

The results for Hochschild homology groups (and hence also for cyclic
homology groups) strongly depend, however, on the particular structure
of the given groupoid. The previous two papers compute these groups
for families of groupoids whose Lie algebroids are rationally
isomorphic to the tangent bundle. This includes families of manifolds
without boundary, families of $b$-pseudodifferential operators.  We
also treated in \cite{BenameurNistor2} the case of families of
pseudodifferential operators on manifolds with corners.

In this paper, we study the Hochschild homology of algebras $\CAA$ of
complete symbols on a foliated manifold $(M,F)$. These algebras can
also be defined using groupoids, although in this paper we choose to
define them directly (see Section \ref{Sec.CSF}). We again obtain a
convergent spectral sequence relating the Hochschild homology groups
of $\CAA$, denoted $\HH_k(\CAA)$, $k = 0,1, \ldots , $ to the Poisson
homology of $(M,F)$ and to the Poisson homology of various other
foliations associated to $(M,F)$. This leads to a complete
determination of the traces of $\CAA$. These traces are usually called
``residue traces.'' See \cite{BenameurNistor1, BenameurNistor2, 
BrylinskiGetzler, Kassel, LauterMoroianu, MelroseNistor, Wodzicki,
WodzickiU} for previous results of this kind.
A motivation for the study of residue traces is our desire to
understand an index theorem of Piazza for pseudodifferential operators
on manifolds with boundary \cite{Piazza}. See also \cite{NT1, NT2}.

For simplicity, we have restricted ourselves here to foliated
manifolds without boundary.  In fact, the first three sections of this
paper are devoted to the study of the Poissson homology of foliated
manifolds and to their relation to longitudinal de Rham cohomology, as
follows. We begin by reviewing some properties of de Rham cohomology
groups for foliations in Section \ref{Sec.deRham}. Then, we discuss in
Section \ref{Sec.Leray} a Gysin long exact sequence for sphere
fibrations of foliated manifolds following \cite{Roger}. The
homogeneous Poisson homology for conic foliated manifolds is defined
and studied in Section \ref{Sec.Can.Hom}.  These homogeneous Poisson
homology groups turn out to be isomorphic to certain de Rham
cohomology groups, see Theorem \ref{hom.can}. The corresponding result
for non homogeneous homologies holds only for the longitudinal Poisson
homology groups, see Definition \ref{eq.HdF} and \cite{VaismanBook}.

In Section \ref{Sec.CSF}, we introduce the algebra $\CAA$ of complete
symbols on a foliated manifold $(M,F)$. We then compute in the last
section, Section \ref{Sec.PSDO}, the $E^2$-term of a spectral
sequence $\EH^r$ that we prove to converge to the Hochschild homology
of longitudinal complete symbols. These computations show that if $p$
and $q$ are, respectively, the dimension and the codimension of the
foliation $(M,F)$, then the groups
$$
	\HH_k(\CAA) = 0\,, \quad \text{if} \quad k > 2p + q.
$$
When the spectral sequence collapses at $E^2$, we get a complete
computation. More precisely, in this case, the Hochschild homology
groups of longitudinal complete symbols are given by (see Corollary
\ref{cor.comp.HH}):
$$
	\HH_k(\maA(M,F)) \simeq \bigoplus_{j=0}^q
	\cohom^{2p+j-k,j}(\S^*F\times \S^1, F_1),
$$
where $\S^*F$ is the longitudinal cosphere bundle and $F_1$ is the
usual foliation on the total space of the bundle $\S^*F\times \S^1$ 
induced by $F$ (with same codimension). 

For the lowest and the highest possibly non-trivial Hochschild
homology groups, our results on the spectral sequence $\EH^r$, when
combined with the Gysin exact sequence mentioned above, show
that
$$
	\HH_0(\maA(M,F)) \simeq \cohom^{p,0}(M,F) \;\;\text{ and }\;\;
	\HH_{2p+q}(\maA(M,F)) \simeq \cohom^{0,q}(M,F), \quad p \ge 2.
$$
(See Theorem \ref{theorem.traces}.)  As a consequence, we obtain a
bijective correspondence between residue traces and holonomy invariant
transverse distributions, as expected.

In the last section, we determine the groups
$\HH_k(\CAA)$ in the following particular case.  Let $M = (\S^1)^n$ be
foliated by the one parameter subgroups $(e^{2\pi \imath \alpha_1 t},
e^{2\pi \imath \alpha_2 t}, \ldots, e^{2\pi \imath \alpha_n t})$. We
assume that the following Diophantine condition is satisfied:\ there
exists $C > 0$ and $N \in \NN$ such that
\begin{equation}\label{eq.diophantine'}
	| m_1 \alpha_1 + m_2 \alpha_2 + \ldots + m_n \alpha_n |^{-1}
	\le C(|m_1| + |m_2| + \ldots |m_n|)^{N},
\end{equation}
for any $m_1, \ldots, m_n \in \ZZ$, not all zero. Then the Hochschild
cohomology groups of this algebra are given by
$$
	\HH^l(\CAA) \cong \Lambda^l \CC^{n+1} \otimes \CC^{\{\pm\}}.
$$
Note that even in this simple example, the assumption of Equation
\eqref{eq.diophantine'} is necessary for this determination to hold. In
general, we need another formulation. This is in sharp contrast with
the behavior of periodic cyclic homology groups.

We use several types of cohomology groups in this paper. The most
important ones are introduced as follows:
\begin{itemize}
\item the longitudinal de Rham cohomology groups $\cohom^{r,s}(M,F)$,
$\cohom^{r,s}(M,F)_l$, and $\cohom_c^{r,s}(M,F)$ are introduced in
Definition \ref{def.bigraded.d.h};
\item the groups $\cohom^k(M,F)$ are introduced using Equation
\eqref{eq.Hk};
\item the definition of the Poisson homology groups
$\cohom_k^\delta(X)$ is recalled in Definition \ref{def.eq2and3nh};
\item the $l$th homogeneous Poisson homology groups
$\cohom^\delta_k(X)_l$ and $\cohom_{r,s}^{\delta_F}(X,\maF)_l$ are
introduced in Definition \ref{def.eq2and3}; and, finally,
\item the groups $\cohom_{k}^{\delta_F}(X)$ are introduced using
Equation \eqref{eq.HdF}.
\end{itemize} 

We assume $M$ to be compact for simplicity. Most of the following
results and constructions work for $M$ non-compact by using cohomology
with compact support.  The proof is the same but notationally more
complicated. In particular, the main computations of Hochschild
homology, Theorems \ref{theorem.comp.hom} and \ref{theorem.traces}
remain true by considering compactly supported cohomology groups.

{\em Acknowledgements.} We thank Robert Lauter, Sergiu Moroianu, Jean
Renault, Claude Roger, and Georges Skandalis for useful discussions.
As we completed our manuscript, we received the preprint \cite{LauterMoroianu2},
which deals with some related questions.

\section{de Rham cohomology for foliations\label{Sec.deRham}}

Throughout this paper, we shall denote by $(M,F)$ a smooth manifold
$M$ of dimension $n$ equipped with a smooth foliation $F$. So $F$ is,
by definition, a smooth, integrable sub-bundle of the tangent bundle
$TM$. The transverse bundle to the foliation $(M,F)$ is the quotient
vector bundle $\nu=TM/F$. We denote by $p$ the dimension of $F$ and by
$q$ the codimension of $F$.  Thus $n = p+q$.

The sections of the longitudinal bundle $F$ will be called
longitudinal vector fields. The sections of the exterior powers
$\Lambda^rF^*$ of the dual vector bundle $F^*$ will be called {\em
longitudinal} differential $r$-forms. The space of longitudinal
differential $r$-forms will be denoted by $\Omega^{r,0}(M,F)$ while
$\Omega^r(M)$ will denote, as customary, the space of differential
$r$-forms on the smooth manifold $M$.  Every longitudinal vector field
on $(M,F)$ is also a vector field on $M$ in the usual sense,
therefore, any differential form on $M$ restricts to a longitudinal
differential form on $(M,F)$. This defines surjections
$$
	\Omega^r(M) \longrightarrow \Omega^{r,0}(M,F).
$$
More generally, a section of the bundle $\Lambda^rF^* \otimes
\Lambda^s\nu^*$ will be called {\em a differential form of bi-degree
$(r,s)$}, or {\em $(r,s)$-differential form}, for short. We denote the
space of $(r,s)$-differential forms on $M$ by $\Omega^{r,s}(M,F)$.

Any choice of a supplementary sub-bundle $H$ to $F$ in $TM$ induces
splittings
\begin{equation}\label{splitting}
	\Theta_H: T^*M  \cong F^* \oplus \nu^* \quad {\text{and}} \quad
	\Omega^d(M) \cong \bigoplus_{r+s=d} \Omega^{r,s}(M,F),
\end{equation}
obtained from the induced embeddings
$$
	\Theta_H = \Theta_H^{r,s} : \Omega^{r,s}(M,F) \longrightarrow
	\Omega^{r+s}(M).
$$
Note that $\nu^*$ identifies canonically with a sub-bundle of $T^*M$
(more precisely, with the annihilator of $F$).  The splitting
\eqref{splitting} endows $\Omega^*(M)$ with a bi-grading so that the de
Rham differential decomposes as a sum of three bi-homogeneous
components
\begin{equation}\label{eq.ddd}
	d=d_F+d_{\perp}+\pa
\end{equation}
where $d_F$ is the $(1,0)$--component, called the longitudinal
differential, $d_{\perp}$ is the $(0,1)$--component and $\pa$ is an
additional map that can be shown to have bi-degree $(-1,2)$
\cite[page 35]{Tondeur}.  Moreover, $d_F$ does not depend on the
choice of the complement $H$ to $F$ in $TM$, as we shall prove
shortly.

In applications, spaces of compactly supported functions are also
needed. Our constructions extend to this case with very few changes.
For simplicity, we shall not consider this case separately.

Let $Z\in \Gamma(\nu)$ be a section of the bundle $\nu$. We shall however
denote by \ $Z_H$ the vector field in $\Gamma(H)$ that corresponds to
$Z$ under the isomorphism $\nu \cong H$. Also, we shall denote by
$\pi_F$ the projection $TM \to F$ with kernel $H$. Let $\theta$ be the
smooth section of $F\otimes \Lambda^2\nu^*$ given for $Y,Z\in
\Gamma(\nu)$ by
$$
	\theta(Y,Z)=\pi_F([Y_H,Z_H]).
$$
Recall that $\pa$ is the contraction by $\theta$, see \cite{Tondeur}
and also \cite[page 267]{ConnesBook}.

The equality $d^2=0$ is then equivalent to
\begin{multline}\label{partial}
	d_F^2=0, \;\; \pa^2=0, \;\; d_{\perp}^2 + \pa d_F + d_F \pa
	=0, \\ d_F d_{\perp}+ d_{\perp}d_F = 0 \;\; \text{ and } \;\;
	\pa d_{\perp}+ d_{\perp}\pa =0.
\end{multline}

Thus, for any $s \in \{0, \ldots, q\}$, we get the complex
\begin{equation}\label{eq.Hrs} 
	0 \to \Omega^{0,s}(M,F) \stackrel{d_F}{\longrightarrow}
	\Omega^{1,s}(M,F) \stackrel{d_F}{\longrightarrow} \ldots
	\stackrel{d_F}{\longrightarrow}\Omega^{p,s}(M,F) \to 0\,,
\end{equation}
called the {\em longitudinal de Rham complex}. If $M$ is endowed with
a free action of $\RR_{+}^*$, we shall denote by
$$
	\Omega^{r,s}(M,F)_l \subset \Omega^{r,s}(M,F)
$$ 
the subspace of forms that are homogeneous of degree $l$ with respect
to the action of $\RR_+^*$. We can assume that the action of $\RR_+^*$
preserves the bundle $H$, and hence that the isomorphism $\Theta_H$ is
invariant with respect to $\Theta_H$.

\begin{lemma}\label{lemma.unique}\  The differential $d_F$ does not
depend on the particular choice of $H$.
\end{lemma}

\begin{proof}\  
%
Let $H$ be the complement of $F$ used to define $d_F$. Denote by
$\pi_{\nu}:TM \to \nu \simeq H$ the quotient projection. For any $Y
\in \Gamma(\nu)$, let $Y^H \in \Gamma(H)$ be the lifting of $Y$ to a
vector field on $M$ such that $\pi_{\nu}(Y^H) = Y$.

Then the differential $d_F$ is explicitly given by
\begin{multline*}
	d_F \omega(X_1,\ldots, X_{r+1};Y_1, \ldots, Y_s) \\ =
	\sum_{j=1}^{r+1} (-1)^{s+j+1} X_j \omega(X_1,\ldots,{\hat
	X}_j, \ldots, X_{r+1};Y_1, \ldots, Y_s) \\ + \sum_{1\leq i < j
	\leq r+1} \omega([X_i,X_j], X_1, \ldots, {\hat X}_i, \ldots,
	{\hat X}_j,\ldots, X_{r+1}; Y_1, \ldots, Y_s) \\ +
	\sum_{i=1}^s \sum_{j = 1}^{r+1} (-1)^{s+j+i}
	\omega(X_1,\ldots,{\hat X}_j, \ldots, X_{r+1};
	\pi_{\nu}[Y^H_i,X_j], Y_1, \ldots,{\hat Y}_i, \ldots, Y_s).
\end{multline*}
where $\omega \in \Omega^{r,s}(M,F)$, $X_1,\ldots, X_{r+1} \in
\Gamma(F)$ and $Y_1, \ldots, Y_s \in \Gamma(\nu)$. Therefore, the only
contribution of the splitting appears in the vector field
$\pi_{\nu}[Y^H_i,X_j]$. But the projection $\pi_{\nu}[Y^H_i,X_j]$
actually does not depend on the particular choice of $H$, because $F$
is integrable. This completes the proof.
\end{proof}

\begin{definition} \label{def.bigraded.d.h}
The $r^{th}$ cohomology space of the longitudinal de Rham complex
\eqref{eq.Hrs} will be denoted by $\cohom^{r,s}(M,F)$. Similarly, we
shall denote by $\cohom^{r,s}(M,F)_l$ the cohomology of the subcomplex
of \eqref{eq.Hrs} consisting of $l-$homogeneous forms and by
$\cohom_c^{r,s}(M,F)$ the cohomology of the subcomplex of
\eqref{eq.Hrs} consisting of compactly supported forms.
\end{definition}

We shall refer to all these groups as the {\em longitudinal de Rham cohomology
groups.}

We shall also need the global longitudinal de Rham complex:
\begin{equation}\label{eq.Hk}
	0 \to \Omega^{0}(M) \stackrel{d_F}{\longrightarrow}
	\Omega^1(M) \stackrel{d_F}{\longrightarrow} \ldots
	\stackrel{d_F}{\longrightarrow}\Omega^n(M) \to 0,
\end{equation}
whose $k^{th}$ cohomology space is denoted $\cohom^k(M,F)$.
So, using  the splitting
\eqref{splitting}, we have:
\begin{equation}
	\cohom^k(M,F) \simeq \bigoplus_{r+s=k} \cohom^{r,s}(M,F).
\end{equation}
The de Rham cohomology spaces of the smooth manifold $M$
will be denoted by $\cohom^k(M)$. 

The homogeneity of $d_F$, $d_\perp$, and $\pa$ gives that
\begin{equation}
	d_F(\omega \wedge \eta) = d_F\omega \wedge \eta + (-1)^r
	\omega \wedge d_F \eta,
\end{equation}
where $\omega \in \Omega^{r,s}(M,F)$ and $\eta \in \Omega^{r',s'}(M,F)$. Since 
$$
	\Theta_H(\Omega^{r,s}(M,F)) \wedge
	\Theta_H(\Omega^{r',s'}(M,F)) \subset
	\Theta_H(\Omega^{r+r',s+s'}(M,F)),
$$ 
 we obtain a product
\begin{equation}\label{eq.prod}
	\cohom^{r,s}(M,F) \otimes \cohom^{r',s'}(M,F) \longrightarrow
	\cohom^{r+r',s+s'}(M,F).
\end{equation}

We shall also need functoriality properties for the groups
$\cohom^{r,s}(M,F)$.

\begin{proposition}\label{prop.funct}\ Let $f : (M,F) \to (M_1,F_1)$
be a $\CI$-map of foliated manifolds such that there exists
complements $H$ and $H_1$ of $F$, respectively $F_1$, with $f_*(H)
\subset H_1$. Then $f$ induces a map
$$
	f^* : \cohom^{r,s}(M_1,F_1) \longrightarrow \cohom^{r,s}(M,F).
$$
\end{proposition}

\begin{proof}\ 
The assumption that $f: (M,F) \to (M_1,F_1)$ is a smooth map of
foliated manifolds implies that $f$ induces a vector  bundle morphism
$f_* : TM \to TM_1$ such that $f_*(F) \subset F_1$. The assumption
that $f_*(H) \subset H_1$ then yields a map $\Gamma(H_1^*) \to \Gamma(H^*)$.  
Together with $\Gamma(F_1^*) \to \Gamma(F^*)$, these two maps give rise to the
map
$$
	f_{r,s}^* : \Omega^{r,s}(M_1,F_1) \longrightarrow
	\Omega^{r,s}(M,F).
$$
Clearly $\Theta^{r,s}_H \circ f_{r,s}^* = f^* \circ \Theta^{r,s}_H$.
Since $f^* \circ d = d \circ f^*$, by checking bi-degrees we see that
$$
	d_F \circ f_{r,s}^* = f_{r+1,s}^* \circ d_{F_1}.
$$
This shows that the maps $f_{r,s}^*$ define a morphism of complexes, and hence
they give rise to a map $f^* : \cohom^{r,s}(M_1,F_1) \to
\cohom^{r,s}(M,F)$, as claimed.
\end{proof}

Functoriality combine with the products \eqref{eq.prod} to define
(external) products
\begin{equation}\label{eq.eprod}
	\cohom^{r,s}(M,F) \otimes \cohom^{r',s'}(M_1,F_1) \to
	\cohom^{r+r', s+s'}(M \times M_1, F \times F_1).
\end{equation}

Note that $\cohom^{0,s}(M,F)$ coincides with the space
$\Omega^s_{bas}(M,F)$ of differential $s$-forms which are basic for
the foliation, i.e. forms $\omega$ such that
$$
	i_Y\omega = 0 \text{ and } L_Y\omega=0, \quad \forall Y \in
	\Gamma(F).
$$
 When restricted to differential forms of bi-degree $(0,*)$, the de
Rham differential $d$ coincides with the sum $d_F + d_{\perp}$ so,
using the equalities \eqref{partial}, we see that $d$ induces a well
defined differential on $\cohom^{0,s}(M,F)$ that coincides with the
differential induced by $d_{\perp}$. Thus the basic complex associated
to $(M,F)$ is given by:
\begin{equation}
	0\to \cohom^{0,0}(M,F) \stackrel{d=d_{\perp}}{\longrightarrow}
	\cohom^{0,1}(M,F)\stackrel{d_{\perp}}{\longrightarrow} \ldots
	\stackrel{d_{\perp}}{\longrightarrow} \cohom^{0,q}(M,F) \to 0.
\end{equation}
The cohomology of this complex will be called, as customary, the {\em
basic de Rham cohomology of the foliated manifold $(M,F)$} and will be
denoted by $\cohom^*_{bas}(M,F)$.

We shall use basic forms to study the behavior of the cohomology
groups $\cohom^{r,s}$ with respect to some fibrations of foliated
manifolds.  To this end, we shall use a Leray spectral sequence with
coefficients in the sheaf of germs of basic forms and the following
well known result of I. Vaisman \cite{VaismanBook}.

\begin{proposition}\label{prop.basic}\ The sequence of sheaves 
$$
	0\to \Omega^{h}_{bas} \longrightarrow \Omega^{0,h}
	\stackrel{d_F}{\longrightarrow}\Omega^{1,h}
	\stackrel{d_F}{\longrightarrow} \ldots
	\stackrel{d_F}{\longrightarrow} \Omega^{p,h} \to 0
$$
is a fine resolution of the sheaf $\Omega^{h}_{bas}$ of basic
$h$-forms.
\end{proposition}

Therefore, the space $\cohom^{r,s}(M,F)$ can be identified with the
$r^{th}$ cohomology space of $M$ with coefficients in the sheaf
$\Omega^{s}_{bas}$.

\begin{corollary}\label{cor.cohom}\ For $0\leq r\leq p$ and 
$0\leq s\leq q$, we have
$$
	\cohom^{r,s}(M,F)\cong \cohom^r(M,\Omega^{s}_{bas}).
$$
\end{corollary}

Similarly, $\cohom_c^{r,s}(M,F)\cong \cohom_c^r(M,\Omega^{s}_{bas}).$

\section{A  Gysin exact sequence\label{Sec.Leray}}

Let now $\pi:E \to M$ be a fiber bundle over the foliated manifold
$(M,F)$. Let $F_E$ be the integrable sub-bundle of the tangent bundle
$TE$ defined by
$$
	F_E := \Ker(p\circ \pi_*)
$$
where $p:TM \to TM/F$. To compute the bi-degree $(r,s)$
cohomology-spaces of the foliated manifold $(E,F_E)$, we can use a
Gysin spectral sequence for the sheaf $\Omega^{s}_{bas}$ over $E$.

The $E_2$ term of this spectral sequence is given by
$E_2^{u,v}=\cohom^u(M,R^v\pi_*(\Omega^{h}_{bas})),$
where the sheaf $R^v\pi_*(\Omega^{h}_{bas})$ is defined by
\begin{equation}
	[R^v\pi_*(\Omega^{h}_{bas})](U)=
	\cohom^v(\pi^{-1}(U),\Omega^{h}_{bas}).
\end{equation}
 Let us now recall the following
result from \cite{Roger}, whose proof we include for the benefit of
the reader.

\begin{proposition}[Roger]\ Let $\pi:E\to M$ be any fibre bundle over $M$
with typical fiber a connected manifold $Y$. If $\mathcal H^v$ denotes
the locally constant presheaf on $M$ defined by ${\mathcal H}^v(U) =
\cohom^v(\pi^{-1}(U))$ and $h$ is arbitrary, but fixed, then there exists a
spectral sequence with
$$
	E_2^{u,v}\cong \cohom^u(M,\Omega^{h}_{bas} \otimes {\mathcal H}^v)
$$
and convergent to $\cohom^{u+v,h}(E,F_E)$.
\end{proposition}

\begin{proof}\ 
Recall that a {\em distinguished open covering} of a foliated manifold
$(M,F)$ is a covering of $M$ by open sets such that the induced
foliation on each of these open sets is a product foliation with
contractable fibers and contractable base. We can always find a
distinguished open covering of the manifold $M$ that also trivializes
the fibration $\pi : E \to M$.  But with respect to any distinguished
open set $U\cong W\times T$, $T$ transversal, such that $\pi^{-1}(U)\cong U\times Y\cong
W\times T\times Y$, we have
$$
	\Omega^{h}_{bas}(W\times T\times Y)\simeq
	\Omega^{h}_{bas}(W\times T)\otimes \cohom^0(Y).
$$
This gives 
	$\cohom^u(\pi^{-1}(U), \Omega^{h}_{bas})\simeq
	\oplus_{u_1+u_2=u} \cohom^{u_1} (W\times T,\Omega^{h}_{bas})
	\otimes \cohom^{u_2}(Y,\R).$
On the other hand, we have:
$\cohom^{u_1,h}(W\times T) \simeq 0 {\mbox{ if }} u_1>0,$
and 
$$
	\cohom^{0,h}(W\times T) \simeq
	\Omega^{h}_{bas}(W\times T)\simeq \Omega^h(T).
$$
We thus obtain
$\cohom^u(\pi^{-1}(U), \Omega^{h}_{bas})\simeq
	\Omega^{h}_{bas}(W\times T) \otimes \cohom^u(Y,\R).$
Our spectral sequence is the spectral sequence associated to the
covering by the open sets $U$ above, and hence
$$
	E_2^{u,v}\simeq \cohom^u(M, \Omega^{h}_{bas} \otimes {\mathcal
	H}^v)
$$
\end{proof}

Now, let $E\stackrel{\pi}{\rightarrow}M$ be an oriented bundle with
fiber of dimension $r$. Denote by $\pi^*$ the pull-back of
differential forms and by $\pi_*$ integration along the fibres of
$E\to M$. If $H$ is a splitting in $(M,F)$ as in \eqref{splitting},
then $\pi^*(H)$ is a splitting for $(E,F_E)$. We fix these splittings
in what follows.

\begin{lemma}\
(i) If $d_{F_E}$ is the longitudinal differential on the foliated
manifold $(E,F_E)$, then
$$
	d_{F_E}\circ \pi^*=\pi^*\circ d_F.
$$

(ii) Similarly, integration along the fibres $\pi_*$ satisfies
$$
	d_F\circ \pi_*=\pi_*\circ d_{F_E}.
$$
\end{lemma}

\begin{proof}\ (i)  follows from Proposition \ref{prop.funct}.

(ii) In the same way, we deduce from the definition that $\pi_*$
is of bi-degree $(-r, 0)$, namely
$$
	\pi_*:\Omega^{k,h}(E,F_E) \to \Omega^{k-r,h}(M,F).
$$
Therefore, from the classical relation $\pi_*\circ d = d\circ \pi_*$
we deduce again by checking bi-degrees that $d_{F}\circ
\pi_*=\pi_*\circ d_{F_E}.$
\end{proof}

\begin{remark}\ Assume that the fibers of $\pi : E \to M$ are
diffeomorphic to the sphere $\S^r$ and that the fibration $E \to M$ is
oriented. The Euler class $e \in \cohom^{r+1}(M)$ is then defined (see
\cite{MilnorStasheff} for details). Moreover, it can be represented by
an element of $\cohom^{r+1,0}(M,F)$ (this actually follows from the
proof of Theorem \ref{Gysin}). Therefore, we have for any
$\alpha\in \Omega^{k,h}(M,F)$ that 
$$
	d_F(\alpha \wedge e) = d_F(\alpha)\wedge e.
$$ 
\end{remark}

In the sequel, in order to make our results more explicit, we shall
need a Gysin exact sequence for the $(k,h)$-cohomology groups. More
precisely, we have

\begin{theorem}\label{Gysin}\ 
Assume that $\pi:E\to M$ is an oriented sphere bundle with fiber
$\S^r$ and denote by $e\in \cohom^{r+1,0}(M)$ the Euler class of this
bundle, then, for any $h$ in the range $0 \le h \le q$, we have the
following Gysin exact sequence
\begin{multline*}
	\ldots\stackrel{\pi^*}{\longrightarrow} \cohom^{k,h}(E,F_E)
	\stackrel{\pi_*}{\longrightarrow} \cohom^{k-r,h}(M,F)
	\stackrel{\wedge e}{\longrightarrow} \cohom^{k+1,h}(M,F)
	\stackrel{\pi^*}{\longrightarrow} \\ \cohom^{k+1,h}(E,F_E)
	\stackrel{\pi_*}{\longrightarrow} \ldots
\end{multline*}
\end{theorem}

\begin{proof}\ 
Because $E$ is an oriented bundle, the presheaf ${\mathcal H}^v$ has
no monodromy (that is, it is constant). Thus, we obtain
$$
	E_2^{u,v}\cong \cohom^{u,h}(M)\otimes \cohom^v(\S^r).
$$ 
Now since $\cohom^v(\S^r)=0$ if $v\not =0$ and $v \not = r$, we get
inclusions $E_{\infty}^{k-r,r} \hookrightarrow E_2^{k-r,r}$,\ $\forall
k\geq 0$. 
Therefore the following sequence is exact
$$
	0\to E_{\infty}^{k-r,r} \hookrightarrow E_2^{k-r,r}
	\stackrel{d_{r+1}}{-\!\!-\!\!-\!\!-\!\!\longrightarrow} E_2^{k+1,0}
	\longrightarrow E_{\infty}^{k+1,0} \to 0.
$$
On the other hand we have an exact sequence 
$$
	0\to E_{\infty}^{k,0}\longrightarrow \cohom^{k,h}(E,F_E)
	\longrightarrow E_{\infty}^{k-r,r} \to 0.
$$
But
$E_2^{k-r,r}\simeq \cohom^{k-r,h}(M,F) {\text{ and }}
	E_2^{k+1,0} \simeq \cohom^{k+1,h}(M,F).$
As in the case of the classical Gysin sequence, the above two short exact
sequences yield a long exact sequence,
\begin{multline*}
	\cdots \longrightarrow \cohom^{k,h}(E,F_E) \longrightarrow
	\cohom^{k-r,h}(M,F) \longrightarrow \cohom^{k+1,h}(M,F)
	\longrightarrow\\ \cohom^{k+1,h}(E,F_E) \longrightarrow
	\cohom^{k+1-r,h}(M,F) \longrightarrow \cdots
\end{multline*}
To end the proof, we must identify the maps involved in this exact
sequence. But this is again similar to the computation for the
classical Gysin sequence.
\end{proof}

From this theorem we obtain the following corollaries.

\begin{corollary}\label{Cor.2}\ We use the notation of Theorem \ref{Gysin}.

(i) The map $\pi^*:\cohom^{k,h}(M,F)\longrightarrow \cohom^{k,h}(E,F_E)$
is an isomorphism for any $r\geq 1$ and $0\leq k\leq r-1$.

(ii)  The map $\pi_*:\cohom^{k,h}(E,F_E) \longrightarrow
\cohom^{k-r,h}(M,F)$ is an isomorphism for any $k \geq p+1$,\
$p = \dim (F)$. 
\end{corollary}

\begin{proof}\ This is a corollary of the longitudinal
Gysin exact sequence proved in Theorem \ref{Gysin}. More precisely,
for $k\leq r-1$ we have $\cohom^{k-r,h}(M,F) = 0.$
Therefore, we get:
\begin{multline*}
	 \ldots \stackrel{\pi_*}{\longrightarrow}
	\cohom^{k-r-1,h}(M,F) =0 \stackrel{\wedge e}{\longrightarrow}
	\cohom^{k,h}(M,F) \stackrel{\pi^*}{\longrightarrow}
	\cohom^{k,h}(E,F_E)\\ \stackrel{\pi_*}{\longrightarrow}
	\cohom^{k-r,h}(M,F) =0 \stackrel{\wedge e}{\longrightarrow}
	\ldots
\end{multline*}

In the same way, if $k\geq p+1$, then $\cohom^{k,h}(M,F)=0$, therefore we get:
\begin{multline*}
	 \ldots \stackrel{\wedge e}{\longrightarrow}
	\cohom^{k,h}(M,F)=0 \stackrel{\pi^*}{\longrightarrow}
	\cohom^{k,h}(E,F_E) \stackrel{\pi_*}{\longrightarrow}
	\cohom^{k-r,h}(M,F)\\ \stackrel{\wedge
	e}{\longrightarrow}\cohom^{k+1,h}(M,F) =0
	\stackrel{\pi^*}{\longrightarrow}\ldots
\end{multline*}
\end{proof}

In particular, for the product $E = M \times \S^{r}$, we get the
following isomorphism that will be used later on.

\begin{corollary}\label{cor.Sr}\
If $E = M \times \S^r$, then
$$
	\cohom^{k,h}(E,F_E) \simeq \cohom^{k,h}(M,F) \oplus
	\cohom^{k-r,h}(M,F),
$$
naturally.
\end{corollary}

\begin{proof}\ The Euler class $e$ vanishes because $E$ is a
product, and hence the Gysin long exact sequence of Theorem
\ref{Gysin} decomposes as a direct sum of short exact sequences
$$
	0 \to \cohom^{k,h}(M,F) \longrightarrow \cohom^{k,h}(E,F_E)
	\longrightarrow \cohom^{k-r,h}(M,F) \to 0.
$$
To complete the proof, it is enough to prove that the above sequence
splits naturally. To this end, let $\omega_r$ be the generator of
$\cohom^r(\S^r)$.  We can pull this class to a cohomology class in
$\cohom^{r,0}(E,F_E)$, denoted $\eta_r$. Then the product with
$\eta_r$ defines the desired natural splitting $\cohom^{k-r,h}(M,F)
\to \cohom^{k,h}(E,F_E)$.
\end{proof}

For $\alpha >0$, the vector bundle $|\Lambda|^{\alpha}(M)$ of
$\alpha$-densities over $M$ is, by definition, the line bundle whose
fiber at a point $x$ is the 1-dimensional complex vector space of maps
$\rho:\Lambda^n(T_xM) \rightarrow \C$ that satisfy
$$
	\forall \lambda\in \R, \, \forall v\in \Lambda^n(T_xM), v \not
	= 0, \quad \rho (\lambda v) = |\lambda|^{\alpha} \rho(v).
$$
This bundle admits nowhere vanishing sections and is, in fact,
trivializable, but not in a canonical way, in general.  Denote by
$\C_M$ the complex orientation bundle of $TM$, then we have
$$
	|\Lambda|^1(M)\cong \Lambda^nT^*M\otimes \C_M.
$$

Let now $E$ be a smooth (real) vector bundle over $M$. The space of
compactly supported smooth sections of $E$ is then naturally endowed
with a structure of a locally convex space. The space of generalized
sections of the vector bundle $E$ is by definition the dual space of
the space of compactly supported smooth sections of the vector bundle
$E^*\otimes |\Lambda^1|(M)$, where $E^*$ is the dual vector bundle of
$E$. Hence a distribution on $M$ can also be viewed as a generalized
1-density.  Some functorial properties of generalized sections are
studied in \cite{GuilleminSternberg}.  In particular, the pull-back of
generalized sections is well defined for fibrations (by integration
along the fibers).

A {\em $k$-current on $M$} is a generalized section of the bundle
$\Lambda^{n-k}(T^*M)\otimes \C_M$. So, a $k$-current on $M$ is, by
definition, a continuous linear form on the space
$$
	C_c^{\infty}(M,\Lambda^{n-k}(TM) \otimes \C_M^* \otimes
	|\Lambda^1|(TM)).
$$
But since,
$$
	|\Lambda^1|(TM) \cong \Lambda^n(T^*M)\otimes \C_M \, {\text{
	and }} \, \Lambda^{n-k}(TM) \otimes \Lambda^n(T^*M) \cong
	\Lambda^k(T^*M),
$$
we get
$$
	\Lambda^{n-k}(TM) \otimes \C_M^* \otimes |\Lambda^1|(TM) \cong
	\Lambda^k(T^*M) \otimes \C_M^* \otimes \C_M.
$$
This shows that any $k$-current $\phi$
defines a linear map
$$
	\phi : \CI_c(M, \Lambda^kT^*M) =: \Omega^k(M) \to \CC.
$$

Denote, as before, by $\nu$ the transverse vector bundle
$\nu=TM/F$. We define a $(k,h)$-current as a generalized section of
the bundle
$$
	\Lambda^{p-k}(F^*)\otimes \Lambda^{q-h}(\nu^*)\otimes \C_M.
$$
We shall denote the space of $(k,h)$-currents by $A_{k,h}(M,F).$

\begin{lemma}\label{lemma.traces}\ 
By choosing a transverse distribution $H$, we can view any
$(k,h)$-current as a continuous linear form on the space of compactly
supported differential $(k,h)$-forms.
\end{lemma}

\begin{proof}\ A $(k,h)$-current on $M$ is by definition a continuous 
linear form on the space
$$
	C_c^{\infty}(M,\Lambda^{p-k}(F) \otimes \Lambda^{q-h}(\nu)
	\otimes \C_M^* \otimes |\Lambda^1|(TM)).
$$
The choice of $H$ fixes an isomorphism $TM \cong F\oplus \nu$ so that
\begin{equation*}
	|\Lambda^1|(TM) \cong \Lambda^n(T^*M)\otimes \C_M \cong
	\Lambda^p(F^*) \otimes \Lambda^q(\nu^*) \otimes \C_M.
\end{equation*}
Using 
\begin{equation*}	
	\Lambda^{p-k}(F)\otimes \Lambda^p(F^*) \cong \Lambda^k(F^*)
	\quad \text{ and } \quad \Lambda^{q-h}(\nu) \otimes
	\Lambda^q(\nu^*) \cong \Lambda^h(\nu^*),
\end{equation*}
we obtain that
$$
	\Lambda^{p-k}(F) \otimes \Lambda^{q-h}(\nu )\otimes \C_M^*
	\otimes |\Lambda^1|(TM) \cong \Lambda^k(F^*)\otimes
	\Lambda^h(\nu^*)\otimes \C_M^*\otimes \C_M.
$$
To finish the proof, we use that the bundle $\C_M^*\otimes \C_M$ is
canonically isomorphic to the trivial line bundle.
\end{proof}

The above lemma shows, in particular, that orientation-twisted
$(p-k,q-h)$-differential forms define a pairing with
$(k,h)$-differential forms.  This is, of course, nothing but the
Poincar\'e map.

For a fixed transverse distribution $H$, we define a longitudinal
differential on the space of $(k,h)$-currents, still denoted $d_F$,
satisfying $d_F^2=0$, which again does not depend on the particular
choice of $H$.  This differential is dual to the one defined above on
smooth differential forms and we get in this way longitudinal
complexes $(A_{*,h}(M,F),d_F)_{0\leq h\leq q}$ of currents:
$$
	0\to A_{p,h}\longrightarrow A_{p-1,h}\longrightarrow \cdots
	\longrightarrow A_{0,h}\to 0.
$$
The cohomology of this complex 
will be denoted $\cohom_{*,h}(M,F)$. So we have a duality map
$\cohom_{k,h}(M,F) \to [\cohom^{k,h}(M,F)]'$, where
$\cohom^{k,h}(M,F)$ is endowed with the quotient topology.

 We include now some remarks that are useful for the reader
interested in relating the above constructions to transverse measures
on foliations.

\begin{definition}\ Let $(M,F)$ be a smooth foliated manifold of dimension $p$ 
and codimension $q=n-p$, as before, and let $\nu=TM/F$ be the
transverse vector bundle.

(i) A {\em transverse current} $C$ on $(M,F)$ is a current of
bi-degree $(p,k)$ for $0\leq k\leq q$, i.e. a generalized section of
the bundle $\Lambda^{q-k}(\nu^*) \otimes \C_M$.

(ii) An {\em invariant} current on $(M,F)$ is a current $C$ on $M$
such that $d_F(C)=0$.

(iii) A current which is transverse and invariant is also called a
{\em basic} current.
\end{definition}

Note that a basic current of type $(p,0)$ is automatically closed in
$M$. Note also that invariant currents are are also sometimes called holonomy invariant currents, 
see \cite{AbouqatebKacimi}. The
simplest example of a transverse current is the Ruelle-Sullivan
current associated with any holonomy invariant transverse measure on
$(M,F)$. Recall that a {\em transverse measure} on $(M,F)$ is a
$\sigma$-finite measure on the disjoint union of submanifolds of $M$
which are everywhere transverse to the foliation. A transverse measure
will be called an {\em invariant transverse measure} if it is
invariant under the action of the holonomy pseudogroup
\cite{Plante}. Given an invariant transverse measure $\mu$, we
canonically associate to $\mu$ an element $C_{\mu}$ of
$C^{-\infty}(M,\Lambda^q\otimes
\C_{\nu})=C^{-\infty}(M,|\Lambda^1|(\nu^*))$ by using partitions of
unity. Therefore, if the foliation is oriented, we have $\CC_M =
\CC_\nu$ and $\mu$ gives rise to a basic current of dimension $q$ that
is closed in $M$.

\section{Canonical homology for foliations\label{Sec.Can.Hom}}

We begin this section by recalling the Koszul-Brylinski complex
\cite{Brylinski} of a foliated Poisson manifold and also some of its
properties that will be needed in the sequel. Let $(M,F)$ be a smooth
foliation with $\dim(M)=n$, $\dim(F)=p$, and $\codim(F)=q$, as
before. We are interested in the manifold $X = F^* \smallsetminus M$, the 
dual of $F$ with the zero section (identified with
$M$) removed. Then $X$ acquires a natural foliation $\maF$ of dimension
$2p$ and codimension $q$. Moreover, $X$ admits an additional
structure, that of a ``foliated Poisson manifold,'' which we proceed
now to define. In the whole section $(X,\maF)$ will then be a foliated
manifold whose leaves have dimension $2p$ and codimension $q$. We
shall insist that $X=F^*$ or some submanifold of $F^*$ when necessary.

\begin{definition}\label{def.Poisson.f}\ 
A {\em foliated Poisson structure} on $(X,\maF)$ is a (foliated)
2-tensor $G \in \Gamma(\Lambda^2 \maF) \subset \Gamma(\Lambda^2 TM)$
over $X$ such that the Schouten-Nijenhuis bracket $[G,G]_{SN}$ is
trivial, see \cite{VaismanBook}.
\end{definition}

A foliated Poisson structure $G$ gives rise to a bilinear form $\{
\cdot , \cdot \}$ on the algebra $\maC^{\infty}(X)$ of smooth maps on
$X$, called {\em the Poisson bracket} and defined by the formula
\begin{equation}
	\{f,g\}=i_G(df\wedge dg),
\end{equation}
where $d$ is the de Rham differential on the smooth manifold $X$ and
$i_G$ is the interior product by the 2-tensor $G$. The condition
$[G,G]_{SN} = 0$ then corresponds to the assumption that $\{ \cdot ,
\cdot \}$ defines a Lie algebra structure on $\CI(X)$. Since for any
$f\in \CI(X)$, the map $g \mapsto \{ f , g \}$ is a derivation of the
commutative ring underlying $\CI(X)$, a foliated Poisson structure on
$X$ endows it with the structure of a Poisson manifold.  Note that the
Hamiltonian vector fields associated with a foliated Poisson structure
are tangent to the leaves of the foliation $(X,\maF)$. The symplectic
leaves of a foliated Poisson manifold $(X,\maF,G)$ are contained in
the foliation $\maF$. When this foliation coincides with the original
foliation $\maF$, we say that the Poisson foliated manifold $(X,\maF,G)$
is a longitudinally symplectic foliated manifold.

A Poisson manifold is a foliated Poisson manifold for any regular
foliation that contains the (singular in general) symplectic
foliation.  Foliated Poisson manifolds are especially
interesting when the symplectic foliation of a given Poisson manifold
can be embedded in a regular foliation of {\em small} dimension. A
regular Poisson manifold $M$, i.e. with a regular symplectic
foliation, is a foliated Poisson manifold for the symplectic foliation
itself, but also for any other foliation that contains the symplectic one. 

An important example for our purposes is that of the cotangent
bundle of any smooth foliation. More precisely, let $(M,F)$ be a
smooth foliated manifold and denote by $\pi_\nu : TM \to TM/F =: \nu$ the quotient
map.  Let $X:= F^*$ be the total space of the longitudinal cotangent
bundle to $(M, F)$ and denote by $\pi:X \to M$ the canonical
projection. The kernel of the composite map $\pi_\nu \circ \pi_*: TX \to
TM/F$ is then an integrable sub-bundle $\maF$ of the tangent bundle
$TX$ to $X$. The leaves of the resulting foliation $\maF$ on $X$ are
exactly the restrictions of the bundle $F^*$ to the leaves of $(M,F)$
and so are symplectic manifolds. By putting together the resulting
symplectic $2$-tensors, we obtain a longitudinally symplectic foliated
manifold $(X,\maF)$.

Let now $(X,\maF,G)$ be a general foliated Poisson manifold. The
Poisson differential $\delta$, is defined as for any Poisson manifold
by the formula \cite{Brylinski}
$$
	\delta := i_G \circ d - d \circ i_G : \Omega^k(X)
	\longrightarrow \Omega^{k-1}(X).
$$

We now recall the definition of Poisson homology of the Poisson
foliation $(X,\maF)$.

\begin{definition}\label{def.eq2and3nh}\
We denote by $\cohom_k^{\delta}(X)$ the {\em Poisson homology of
$X$}, defined by
\begin{equation*}
	\cohom_k^{\delta}(X) := {\frac{\Ker(\delta: \Omega^k(X) \to
	\Omega^{k-1}(X))}{\delta(\Omega^{k+1}(X))}}.
\end{equation*}
\end{definition}

Assume that we have fixed a splitting $\Theta_H$ as in
\eqref{splitting} for the foliated manifold $(X,\maF)$. This, in turn, 
fixes isomorphisms $\Omega^k(X) \simeq \bigoplus_{r+s=k} \Omega^{r,s}(X,\maF)$.

\begin{lemma}\label{delta}\cite{VaismanBook}\
Let $(X,\maF,G)$ be a foliated Poisson manifold, then the
Koszul-Brylinski operator $\delta$, has a canonical decomposition into
two bi-homogeneous operators
$$
	\delta = \delta_{\maF}+ \delta_{-2,1},
$$
where $\delta_{\maF}=[i_G,d_{\maF}]$ is a component of bi-degree
$(-1,0)$ with respect to the splitting, called the longitudinal
Poisson differential, and $\delta_{-2,1}$ is an extra term with
bi-degree $(-2,1)$ with respect to the bi-grading. Furthermore, if $d =
d_{\maF} + d_\perp + \pa$, as in Equation \eqref{eq.ddd} with $(X,\maF)$
in place of $(M,F)$, we have
$$
	\delta_{-2,1} = [i_G,d_{\perp}]\,, \;\; \delta_{\maF}^2 = 0\,, \quad
	\delta_{-2,1}^2 = 0\,, \quad \text{ and } \quad \delta_{\maF}\delta_{-2,1} +
	\delta_{-2,1} \delta_{\maF} = 0.
$$
\end{lemma}

\begin{proof}\ 
Let $H$ be a supplementary sub-bundle to $\maF$ in $TX$ and
$d=d_{\maF} + d_{\perp} + \pa$ the corresponding decomposition of the
de Rham operator $d$ as recalled in Section~\ref{Sec.deRham}. Let us
show that $[i_G,\pa]=0$. Let $X\in \Gamma(\maF)$ be a longitudinal
vector field. Then, for any $X_1,...,X_{k-1}\in \Gamma(\maF)$, for any
$Y_1,...,Y_{h+2}\in \Gamma(H)\simeq \Gamma(\nu)$, and for any $\omega\in
\Omega^{k+1,h}(X,\maF)$, we have:
\begin{multline*}
	\pa(i_X\omega)(X_1,...,X_{k-1};Y_1,...,Y_{h+2})=\sum_{1\leq
	j<i\leq h+2} \\ (-1)^{i+j+h}
	\omega(X,\pi_{\maF}[Y_j,Y_i],X_1,...,X_{k-1}; Y_1,...,{\hat
	Y}^j,...,{\hat Y}^i,...,Y_{h+2}).
\end{multline*}
On the other hand:
\begin{multline*}
	i_X(\pa\omega)(X_1,...,X_{k-1};Y_1,...,Y_{h+2})=\sum_{1\leq
	j<i\leq h+2}\\ (-1)^{i+j+h}
	\omega(\pi_{\maF}[Y_j,Y_i],X,X_1,...,X_{k-1}; Y_1,...,{\hat
	Y}^j,...,{\hat Y}^i,...,Y_{h+2}).
\end{multline*}
Thus we deduce that $i_X\circ \pa + \pa \circ i_X=0$ and hence
$$
	[i_{X\wedge Y},\pa]= i_Y (i_X \pa + \pa i_X) - (i_Y \pa + \pa
	i_Y) i_X = 0,
$$
for any $(X,Y)\in \Gamma(\maF)$.  Therefore $[i_A, \pa]=0,$ for all
$A\in \Gamma(\Lambda^2\maF)$.  We finish the proof by setting
$\delta_{-2,1}=[i_G,d_{\perp}]$.  Finally the identity $\delta^2=0$
\cite{Brylinski} gives the claimed equalities by direct inspection of
the bi-degrees.  See also \cite[Proposition 4.13]{VaismanBook}.
\end{proof}

\begin{remark}\
The contraction by $G$ has bi-degree $(-2,0)$ and satisfies the
relation
$$
	i_G(\omega_1\wedge \omega_2)=i_G(\omega_1)\wedge \omega_2,
	\quad \forall \omega_1\in \Gamma(\Lambda^*T^*X) \text{ and }
	\forall \omega_2\in \Gamma(\Lambda^*H^*).
$$
\end{remark}

Assume for the rest of this section that $(X,\maF)$ is a
longitudinally symplectic foliation with $\dim(\maF)=2p$ and
codim$(\maF)=q$.  

If $U$ is a distinguished chart for the foliation $(X,\maF)$, then
$\delta_{\maF}$ restricts to $U$ and induces a well defined
differential on the sheaf of germs of smooth longitudinal differential
forms. The action of $\delta_{\maF}$ on typical longitudinal forms is
similar to the classical one.  More precisely:

\begin{proposition}\label{prop.can.diff}\
Let $(X, \maF,G)$ be a general foliated Poisson manifold.  Then the
action of $\delta_{\maF}$ on typical longitudinal forms is given by
\begin{multline*}
	\delta_{\maF}(f_0 d_{\maF} f_1 \ldots d_{\maF} f_k) =
   	\sum_{1\leq j\leq k} (-1)^{j+1} \{f_0,f_j\}d_{\maF} f_1 \ldots
   	{\widehat{d_{\maF} f_j}} \ldots d_{\maF} f_k \\ + \sum_{1\leq
   	i<j\leq k} (-1)^{i+j}f_0d_{\maF} \{f_i,f_j\} d_{\maF} f_1
   	\ldots {\widehat{d_{\maF} f_i}} \ldots {\widehat{d_{\maF}
   	f_j}} \ldots d_{\maF} f_k,
\end{multline*}
for all $f_0, \ldots , f_k \in \CI(X)$.
\end{proposition}

\begin{proof}\ The computations carried out in  
\cite[page 96]{Brylinski} imply our proposition. Recall that we have
\begin{multline*}
	\delta(f_0 d f_1 d f_2 \ldots d f_k) = \sum_{1\leq j\leq k}
   	(-1)^{j+1} \{f_0,f_j\}d f_1 \ldots {\widehat{d f_j}} \ldots d
   	f_k \\ + \sum_{1\leq i<j\leq k} (-1)^{i+j}f_0d \{f_i,f_j\} d
   	f_1 \ldots {\widehat{d f_i}} \ldots {\widehat{d f_j}} \ldots d
   	f_k,
\end{multline*}
for all $f_0, f_1,\ldots, f_k \in \CI(X)$. Hence taking the $(k-1,0)$
component of each side of the above equality gives exactly the allowed
formula.
\end{proof}

Let us also mention, for completeness, the following result.

\begin{proposition}\
(i) If $\omega\in \Gamma(\Lambda^k\maF^*)$ and $\omega'\in
\Gamma(\Lambda^{k'}H^*)$,  then we have
$$
	\delta_{-2,1}(\omega\wedge \omega')=
	\delta_{-2,1}(\omega)\wedge \omega'.
$$

(ii) For any $X_1,..., X_r\in \Gamma(\maF)$, for any $Z\in \Gamma(H)$,
and for any $\omega\in \Omega^{r+2,0}(X,\maF)$,
$$
	\delta_{-2,1}\omega(X_1,...,X_r,Z)=i_{[G,Z]_{SN}}\omega
	(X_1,...,X_r).
$$
\end{proposition}

\begin{proof}\
(i) For any $\omega\in \Gamma(\Lambda^k\maF^*)$ and any $\omega'\in
\Gamma(\Lambda^{k'}H^*)$, we have
$$
	d_{\perp}(\omega\wedge \omega')=d_{\perp}(\omega)\wedge
	\omega'+(-1)^k\omega \wedge d_{\perp}(\omega').
$$
Therefore 
$$
	[i_G,d_{\perp}](\omega\wedge \omega')=i_G(d_{\perp}\omega
	\wedge \omega')-d_{\perp}(i_G(\omega))\wedge
	\omega'=[i_G,d_{\perp}] (\omega)\wedge \omega'.
$$
(ii) Let $Y_1, Y_2, X_1, \ldots, X_r \in \Gamma(\maF)$, $Z \in \Gamma(H)$, and 
$\omega \in \Omega^{r+2,0}(X,\maF)$ be arbitrary. Using a simple computation, we obtain
\begin{multline*}
	([i_{Y_1\wedge Y_2}, d_{\perp}]\omega) (X_1, \ldots, X_r, Z)
	=\\ \omega(\pi_{\maF}[Z,Y_1], Y_2, X_1,\ldots, X_r)-
	\omega(\pi_{\maF}[Z,Y_2], Y_1, X_1,\ldots, X_r),
\end{multline*}
where $\pi_{\maF}$ is the projection onto $\maF$ along $H$.
Therefore, we get:
$$
	[i_{Y_1\wedge Y_2}, d_{\perp}] = i_{\pi_{\maF}[Z,Y_1]\wedge
	Y_2 - \pi_{\maF}[Z,Y_2]\wedge Y_1}.
$$
By direct inspection from the definition of the Schouten-Nijenhuis
bracket, we deduce that
$$
	[i_{Y_1\wedge Y_2}, d_{\perp}] = i_{\pi_{\Lambda^2 \maF}([Y_1
	\wedge Y_2,Z]_{SN})},
$$
where $\pi_{\Lambda^2 \maF}$ is the projection onto the space of
longitudinal $(2,0)$-vectors.  But since $\omega \in \Omega^{r+2,0}(X,
\maF)$, this completes the proof.
\end{proof}

We continue to assume for the rest of this section that $(X,\maF)$ is a
longitudinally symplectic foliation with $\dim(\maF)=2p$ and
codim$(\maF)=q$.  For any leaf $L$ of the foliation $\maF$ of $X$, let
$\omega_L$ be the symplectic two form of $L$. Then there exists
longitudinal 2-forms on $X$ that restrict on each leaf $L$ to
$\omega_L$.  If we use the splitting \eqref{splitting} then we can
choose in a unique way a differential 2-form $\omega\in
\Omega^{2,0}(X, \maF)$ that restricts to $\omega_L$ on each leaf
$L$. The form $\omega$ will be called the longitudinal symplectic form
of $(X, \maF)$. It depends on the splitting \eqref{splitting}.

Using the longitudinal symplectic form $\omega$ we can recover the
longitudinal volume form associated with the symplectic orientation by
setting:
\begin{equation}
	vol_{\maF}(X) := \frac{1}{p!}\omega^p\,.
\end{equation}
We then define the longitudinal symplectic Hodge operator $*_{\maF}:
\Omega^{r,0}(X, \maF) \to \Omega^{2p-r,0}(X,\maF)$ by the equality:
\begin{equation}
	\beta \wedge (*_{\maF} \alpha) =
	(\beta,\alpha)_{\omega}.vol_{\maF}(X), \quad \forall
	\alpha,\beta \in \Omega^{r,0}(X, \maF),
\end{equation}
where $(\,\cdot\,,\cdot\,)_{\omega}$ is the bilinear form induced by
the symplectic form on longitudinal differential forms.

\begin{remark}\
For any $f\in \CI(X)$ we have by the definition of $*_{\maF}$:
\begin{equation}\label{*}
	*_{\maF} (f \alpha) = f *_{\maF}\alpha.
\end{equation}
\end{remark}

Recall now (Definition \ref{def.bigraded.d.h}) that $\cohom^{r,s}(X,
\maF)$ denotes the $r^{th}$ cohomology group of the longitudinal
complex
\begin{equation}\label{eq.HdRF}
	0\to \Omega^{0,s}(X, \maF)
	\stackrel{d_{\maF}}{\longrightarrow} \Omega^{1,s}(X, \maF)
	\stackrel{d_{\maF}}{\longrightarrow} \ldots
	\stackrel{d_{\maF}}{\longrightarrow} \Omega^{2p,s}(X, \maF)
	\to 0,
\end{equation}
and $\cohom^{\delta_{\maF}}_{r,s}(X, \maF)$ is the longitudinal
Poisson homology of $(X,\maF)$, that is, the $r^{th}$-cohomology
group of the complex
\begin{equation}\label{eq.HdF}
	0 \to \Omega^{2p,s}(X, \maF)
	\stackrel{\delta_{\maF}}{\longrightarrow} \Omega^{2p-1,s}(X,
	\maF) \stackrel{\delta_{\maF}}{\longrightarrow} \ldots
	\stackrel{\delta_{\maF}}{\longrightarrow} \Omega^{0,s}(X,
	\maF) \to 0.
\end{equation}
The cohomology $\cohom_*^{\delta_{\maF}}(X, \maF)$ of the global complex
$(\Omega^k(X))_{0\leq k\leq 2p+q}$ with respect to the operator
$\delta_{\maF}$ is hence given by
$$
	\cohom_k^{\delta_{\maF}}(X, \maF)\simeq \bigoplus_{r+s=k}
	\cohom^{\delta_{\maF}}_{r,s}(X, \maF).
$$
We hope the reader will be able to easily tell apart all these cohomology
groups and distinguish for instance the groups $\cohom_*^{\delta_{\maF}}(X,
\maF)$ from the Poisson homology groups of $X$ that are denoted 
$\cohom_*^{\delta}(X)$, (see the end of the introduction for a list of 
references to the definitions of the main cohomology groups).

When the foliation $(X,\maF)$ is longitudinally symplectic, the
longitudinal symplectic Hodge operator extends to a well defined
operator, still denoted $*_{\maF}$,
$$
	*_{\maF} : \Omega^{r,s}(X, \maF) \to \Omega^{2p-r,s}(X, \maF),
	\quad \forall s\in \{0,...,q\}.
$$
defined by
$$
	*_{\maF}(\alpha \wedge \beta) := *_{\maF}(\alpha) \wedge
         \beta,
$$
for any $\alpha\in \Omega^{r,0}(X, \maF)$ and any $\beta\in
\Omega^{0,s}(X, \maF).$ This is a consequence of the relation
\eqref{*} and the splitting \eqref{splitting}.  We then see that
$*_{\maF}^2 = 1$.

We point out that the longitudinal Poisson differential
$\delta_{\maF}$ also satisfies a similar relation, namely
$$
	\delta_{\maF}(\alpha \wedge \beta) := \delta_{\maF}(\alpha)
	\wedge \beta,
$$
which follows from the formula given for $\delta_{\maF}$ in
Proposition \ref{prop.can.diff} using the same method as in
\cite{BenameurNistor2}.

\begin{proposition}[Vaisman]\label{prop.ext.can}\
Let $(X, \maF)$ be a longitudinally symplectic foliated manifold with
leaves of dimension $2p$.

(1) We have $(-1)^{r+1}*_{\maF} \circ d_{\maF}\circ *_{\maF} =
\delta_{\maF},$ on $\Omega^{r,s}(X, \maF)$.

(2) The cohomology of $X$ with respect to $\delta_{\maF}$ is given by
$$
	\cohom^{\delta_{\maF}}_{r,s}(X, \maF) \cong \cohom^{2p-r,s}(X,
	\maF),
$$
and hence $\cohom_k^{\delta_{\maF}}(X, \maF)\cong \oplus _{0\leq j\leq
k} \cohom^{2p-j,k-j}(X, \maF).$
\end{proposition}

\begin{proof}\ The proof of (1) is in 
\cite[page 80]{VaismanBook}.  It can also be derived easily from the
properties listed above.

(2) We have $\cohom^{\delta_{\maF}}_{k}(X, \maF) \simeq \bigoplus_{k =
r + s}\cohom_{r,s}^{\delta_{\maF}}(X, \maF)$.  The above result (2),
shows that $*_{\maF}$ induces an isomorphism
$\cohom^{\delta_{\maF}}_{r,s}(X, \maF) \simeq \cohom^{2p-r,s}(X,
\maF)$, extending the case $s = 0$. This proves that
$\cohom^{\delta_{\maF}}_k(X,\maF) \simeq \oplus_{j=0}^k
\cohom^{2p-j,k-j}(X, \maF)$, as claimed. See \cite{VaismanBook} again.
\end{proof}

\section{Conic foliations and their cohomology\label{Sec.Conic}}

We now introduce the action of $\RR_+^*$ into the picture.

\begin{definition}\label{conic}\
Let $(X,\maF,G)$ be a longitudinally symplectic foliation.  The triple
$(X,\maF,G)$ will be called a {\em longitudinally symplectic conic
foliation} if there exists a free smooth action of the group $\R^*_+$
on $X$ by leaf-preserving diffeomorphisms $(\alpha_t)_{t>0}$ such that
$(\alpha_t)_*(G)=G/t$.
\end{definition}

This definition means that each leaf is a conic symplectic manifold in
the sense of \cite{BrylinskiGetzler} and that the global action
$\alpha$ is smooth. For $l\in \Z$, recall that a differential form
$\omega\in \Omega^k(X)$ is $l$-homogeneous if
\begin{equation}
	(\alpha_t)^*(\omega)=t^l\omega, \quad  \forall t>0.
\end{equation}
We shall denote as before by $\Omega^k(X)_l$ the space of
$l$-homogeneous differential $k$-forms on $X$. From the definition of
a longitudinally symplectic conic foliation, we deduce that the
longitudinally symplectic form corresponding to the bivector $G$
belongs to $\Omega^2(X)_1$.

Since the action of $\R^*_+$ on the longitudinally symplectic conic
foliation $(X,\maF)$ is free, we can choose the complement $H$ to $F$
in $TX$ to be $\R^*_+$ invariant. The bi-grading on forms is also
$\R^*_+$-equivariant and we shall denote, as before, by
$\Omega^{r,s}(X,\maF)_l$ the smooth $l$-homogeneous sections of
$\Lambda^r\maF^*\otimes \Lambda^s\nu^*$.

Let now $(M,F)$ be a smooth foliated manifold and take $X = F^*
\smallsetminus M$.  We are interested in the foliated manifold $(X,
\maF)$, where $\maF$ is the foliation defined on the total space $F^*$
of the longitudinal cotangent bundle to $(M,F)$ as before and then
restricted to $X=F^*\smallsetminus M$.  The radial
action of $\R^*_+$ allows us to consider $l$-homogeneous forms
$\Omega^k(X)_l$ and $\Omega^{r,s}(X,\maF)_l$. As we have already
observed, the foliated manifold $(X,\maF)$ is then longitudinally
symplectic. The Poisson differential $\delta$ associated with the
Poisson structure of $X$ sends $\Omega^k(X)_l$ to
$\Omega^{k-1}(X)_{l-1}$. The same holds for the operators $\delta_{\maF}$
and $\delta_{-2,1}$ defined in the previous section.

\begin{definition}\label{def.eq2and3}\
We denote by $\cohom_k^{\delta}(X)_l$ the {\em $l$-homogeneous Poisson
homology of $X$}, defined by
\begin{equation*}
	\cohom_k^{\delta}(X)_l := {\frac{\Ker(\delta:\Omega^k(X)_l \to
	\Omega^{k-1}(X)_{l-1})}{\delta(\Omega^{k+1}(X)_{l+1})}}.
\end{equation*}
In the same way, using again the splitting \eqref{splitting}, we set
\begin{multline*}
	\cohom_{r,s}^{\delta_{\maF}}(X,\maF)_l :=
	{\frac{\Ker(\delta_{\maF}:\Omega^{r,s}(X,\maF)_l \to
	\Omega^{r-1,s}(X,\maF)_{l-1})}{\delta_{\maF}
	(\Omega^{r+1,s}(X,\maF)_{l+1})}} \quad \text{ and }\\
	\cohom_k^{\delta_{\maF}}(X,\maF)_l :=
	{\frac{\Ker(\delta_{\maF}:\Omega^k(X)_l \to
	\Omega^{k-1}(X)_{l-1})}{\delta_{\maF} (\Omega^{k+1}(X)_{l+1})}}.
\end{multline*}
\end{definition}

The homogeneous Poisson complex $(\Omega^*(X)_*,\delta)$ splits into a
direct sum of finite homogeneous subcomplexes $(\maP^k)_{k\in \Z}$
defined by:
\begin{equation}
	\maP^k: 0\to \maP^k_{2p+q-k} \stackrel{\delta}{\longrightarrow}
	\maP^k_{2p+q-k-1} \stackrel{\delta}{\longrightarrow} \ldots
	\stackrel{\delta}{\longrightarrow}\maP^k_{-k} \to 0,
\end{equation}
where $\maP^k_l := \Omega^{k+l}(X)_l$. Therefore we have:
$$
	\cohom^{\delta}_{k+l}(X)_l = {\frac{\Ker(\delta:\maP^k_l \to
	\maP^k_{l-1})}{\delta(\maP^k_{l+1})}}.
$$
If we define in the same way $\maP^{r,s}_l := \Omega^{r+l,s}(X,\maF)_l$,
then we get a further splitting:
$$
	\maP^k_l \simeq \bigoplus_{r+s=k} \maP^{r,s}_l.
$$
With respect to this splitting, the differential $\delta_{\maF}$ preserves
$\maP^{r,s} := \oplus_{l\in \Z}\maP^{r,s}_l$ and sends $\maP^{r,s}_l$
to $\maP^{r,s}_{l-1}$. Thus, to compute the homogeneous homology of
the longitudinal Poisson differential $\delta_{\maF}$, we can restrict
ourselves to $\maP^{r,s}$. Note though that the extra differential
$\delta_{-2,1}$ does not preserve $\maP^{r,s}$ and sends it to
$\maP^{r-1,s+1}$.

Our next result is that, in order to compute homogeneous Poisson
homology, we can get rid of the extra term $\delta_{-2,1}$. 

\begin{proposition}\label{prop.d=dF}\ 
Let $(X,\maF)$ be a longitudinally symplectic conic foliation. Then
$$
	\cohom_{k+l}^{\delta}(X)_l \cong
	\cohom_{k+l}^{\delta_{\maF}}(X,\maF)_l.
$$ 
\end{proposition}

\begin{proof}\ Recall that we have:
$$
	\delta = \delta_{\maF} + \delta_{-2,1} \quad \text{ and } \quad 
	\delta_{\maF}\delta_{-2,1} + \delta_{-2,1} \delta_{\maF} = 0.
$$
Thus, for any fixed $k$, we use the decomposition $\maP^k \cong
\bigoplus_{i+j=k}\maP^{i,j}$ into a finite double complex.  We set
for any fixed $k\in \Z$,
$$
	K^{j,l} := \maP^{k-j,j}_{l-j},
$$
so that 
$$
	\delta_{\maF}:K^{j,l} \longrightarrow K^{j,l-1} \text{ and }
	\delta_{-2,1}: K^{j,l} \longrightarrow K^{j+1,l}.
$$
To compute the homogeneous $\delta$-homology of $X$, we use that the
complex splits into the subcomplexes $(\maP^k,\delta)$. Therefore, we
can fix the integer $k\in \Z$ and define a filtration of the above
bicomplex $K^{j,l}$ by setting
$$
	F_{h} := \bigoplus_{l\in \Z, j\leq h} K^{j,l}.
$$
This yields a spectral sequence $(E^r)_{r\geq 1}$ which converges to
the $\delta$-homology because it comes from a filtration that is
bounded both below and above.  The $E^1$ term of this spectral
sequence is computed by a de Rham cohomology group
$$
	E^1_{u,v} \simeq \cohom^{2p-v-k+u,u}(X,\maF)_{p-k+u},
$$
the isomorphism being implemented by the leafwise symplectic duality
operator $*_{\maF}$.  We now observe that the homogeneous longitudinal
de Rham cohomology space $\cohom^{2p-v-k+u,u}(X,\maF)_{p-k+u}$ is
trivial unless $u=k-p$, by the homotopy invariance of de Rham
cohomology. Therefore, we get
$$
	E^1_{u,v}=0 \text{ if } v\not = -k-p.
$$
Hence for any $r\geq 1$, we see that $d^r=0$ and the spectral sequence
collapses at $E^1$. The proof is thus complete since the spectral
sequences considered are convergent.
\end{proof}

\begin{corollary}\label{prop.ext.can.con}\
Let $(X,\maF)$ be a longitudinally symplectic conic foliation with
leaves of dimension $2p$. Then $\cohom^{\delta_{\maF}}_{r,s}(X,\maF)_l
\cong \cohom^{2p-r,s}(X,\maF)_{l+p-r}$, and hence
$$
	\cohom_k^{\delta_{\maF}}(X,\maF)_l \cong \bigoplus _{0\leq j\leq k}
	\cohom^{2p-j,k-j}(X,\maF)_{l+p-j}.
$$
\end{corollary}

\begin{proof}\  This is a consequence of Proposition
\ref{prop.ext.can}. Note that we have:
$$
	\cohom^{\delta_{\maF}}_k(X,\maF)_l \simeq \bigoplus_{r+s=k}
	\cohom^{\delta_{\maF}}_{r,s}(X,\maF)_l.
$$
But by definition of the operator $*_{\maF}$, we see that it sends
$l$-homogeneous forms of bi-degree $(r,s)$ to $(l+p-r)$-homogeneous
forms of bi-degree $(2p-r,s)$.
\end{proof}

In the case of trivial foliations by 2-planes, that is, when
$X=\R^{2}\times \R^q$ foliated by the symplectic planes $\R^{2}\times
\{pt\}$, if we denote by $(x,\xi)$ the symplectic coordinates along
the leaves and $(y_1,...,y_q)$ the transverse coordinates, we have the
following easy generalizations of some equations in
\cite{Brylinski}. Namely, for any $f,g,h \in \CI(X):$
\begin{equation*}
\begin{gathered}
	*_{\maF}(f)=f dx d\xi, \quad *_{\maF}(fdx+gd\xi)=-(fdx+gd\xi),\\
	*_{\maF}(hdy_i)=hdx d\xi dy_i \quad \text{ and } \quad *_{\maF}(f dx
	*d\xi)=f.
\end{gathered}
\end{equation*}
In the same way we have
\begin{equation*}
\begin{gathered}
	*_{\maF}(f dx dy_{i_1} dy_{i_2}\ldots dy_{i_k})= -(f dx dy_{i_1}
	dy_{i_2}\ldots dy_{i_k}),\\ *_{\maF} (f d\xi dy_{i_1} dy_{i_
	2}\ldots dy_{i_k})= -(f d\xi dy_{i_1} dy_{i_2}\ldots
	dy_{i_k}), \\ *_{\maF}(f dy_{i_1} dy_{i_2}\ldots dy_{i_k})=f dx
	d\xi dy_{i_1} dy_{i_2}\ldots dy_{i_k}\\ \text{ and } *_{\maF}(f dx
	d\xi dy_{i_1} dy_{i_2}\ldots dy_{i_k})=f dy_{i_1}
	dy_{i_2}\ldots dy_{i_k}.
\end{gathered}
\end{equation*}
On the other hand,
$$
	\{f,g\}={\frac{\partial f}{\partial \xi}}{\frac {\partial
	g}{\partial x}}-{\frac{\partial f}{\partial x}}{\frac{\partial
	g}{\partial \xi}}.
$$
Hence Propositions \ref{prop.ext.can} and \ref{prop.ext.can.con} can
also be proved by reducing to the above trivial case, as in
\cite{Brylinski}, for example.

\begin{remark}\ 
The operator $\delta_{-2,1}$ gives rise to a new differential
\begin{equation}
	d_{2,1}=(-1)^{r+1}*_{\maF}\circ \delta_{-2,1}\circ *_{\maF}
\end{equation}
on $\Omega^{r,s}(X,\maF)$ whose bi-degree is $(2,1)$ and which
satisfies
\begin{equation}
	d_{2,1}^2 = 0, \quad d_{2,1} d_{\maF} + d_{\maF} d_{2,1} = 0.
\end{equation}
\end{remark}

We are now in position to compute the homogeneous Poisson homology
of a longitudinally symplectic conic foliation.  The homogeneous
Poisson homology spaces $\cohom^{\delta}_{k}(X)_l$ were defined in
Definition \ref{def.eq2and3}.

\begin{theorem}\label{hom.can}\ 
Let $(X, \maF)$ be a longitudinally symplectic conic foliation.  We
denote $\dim (\maF) = 2p$ and $\codim (\maF) = q$, as before. Then
$$
	\cohom^{\delta}_{k}(X)_l \cong \cohom^{p-l,k-l-p}(X,\maF)_0
$$
for $0\leq k \leq 2p+q$ and $|l|\leq p$. For the other values of $k$
and $l$ we have $\cohom^{\delta}_{k}(X)_l = 0$.
\end{theorem}

\begin{proof}\ 
We use Proposition \ref{prop.d=dF} to conclude that
$\cohom_{k+l}^{\delta} (X)_l \simeq \cohom_{k+l}^{\delta_{\maF}}
(X,\maF)_l.$ By Proposition \ref{prop.ext.can.con},
$$
	\cohom_k^{\delta_{\maF}}(X,\maF)_l \cong \oplus_{0\leq j\leq k}
	\cohom^{2p-j,k-j}(X,\maF)_{l+p-j}.
$$
By the homotopy invariance of de Rham cohomology, only the groups for
which $l + p - j = 0$, are non-zero. The result is obtained then by
substituting $j = p +l$. The other groups vanish for dimension
reasons.
\end{proof}

Let now $X = F^*\smallsetminus M$ be the dual of the foliation $F$ of
$M$ with the zero section removed and with the induced structure of a
longitudinal symplectic conic foliation. Let $\S^*F = X/\RR_{+}^*$ be
the cosphere bundle of $F$ and $F_1$ be the induced foliation on
$\S^*F \times \S^1$ with leaves of dimension $2p$ (each copy of $\S^1$
is completely contained in a leaf). The above theorem then gives the
following result.

\begin{corollary}\label{cor.hom.can}\ 
%
Let $(M,F)$ be a foliated manifold and $X = F^*\smallsetminus M$,
$\S^*F = X/\RR_{+}^*$, and $F_1$ be as in the paragraph above. We
denote $p = \dim (F_1)/2 = \dim(F)$ and $q = \codim(F_1) = \codim(F)$. Then
\begin{multline*}
	\cohom^{\delta}_{k}(X)_l \cong \cohom^{p-l,k-l-p}(\S^*F \times
	\S^1,F_1) \\ \cong \cohom^{p-l,k-l-p}(\S^*F, \S^*F \cap F_1)
	\oplus \cohom^{p-l-1,k-l-p}(\S^*F, \S^*F \cap F_1)
\end{multline*}
for $0\leq k \leq 2p+q$ and $|l|\leq p$. For the other values of $k$
and $l$, we have $\cohom^{\delta}_{k}(X)_l = 0$.
\end{corollary}

\begin{proof}\ The proof is exactly as in the case when the foliation 
$F$ is trivial (with just one leaf) \cite{BenameurNistor1,
BenameurNistor2, BrylinskiGetzler, MelroseNistor}. The crucial
ingredient of the proof is to choose a function $r \not = 0$
homogeneous of degree one. Then identify $r^{-1}dr$ with $\omega_1$,
the generator of $H^1(\S^1)$. This leads to the isomorphism
$\cohom^{p-l,k-l-p}(X,\maF)_0 \cong \cohom^{p-l,k-l-p}(\S^*F \times
\S^1,F_1)$.  The second isomorphism follows from Corollary
\ref{cor.Sr}.
\end{proof}

We are ready now to handle an explicit example.

\begin{example}\label{ex.torus}\ 
Let us consider $M = (\S^1)^n$, foliated by the one parameter
subgroups $(e^{2\pi \imath \alpha_1 t}, e^{2\pi \imath \alpha_2 t},
\ldots, e^{2\pi \imath \alpha_n t})$, not all of $\alpha_i$'s equal to
zero. Thus $p = 1$ and $q = n-1$. Then
$$
	F = M \times \RR, \quad X = M \times (\RR\smallsetminus \{0\}),
	\quad \S^*F \times \S^1 \simeq M \times \{\pm \} \times \S^1,
$$
with leaves $L \times \{\epsilon\} \times \S^1$, where $L \subset M$
is a leaf of $M$ and $\epsilon = +$ or $\epsilon = -$. Then the second
isomorphism in Corollary \ref{cor.hom.can} gives
$$
	\cohom^{k,h}(\S^*F \times \S^1, F_1) \cong \bigr(
	\cohom^{k,h}(M,F) \oplus \cohom^{k-1,h}(M,F) \bigl ) \otimes
	\CC^{ \{\pm \} }.
$$
($\CC^{\{\pm\}}$ is the complex vector space with basis $+$ and $-$.)
To obtain more precise results (which happen to also be finite
dimensional spaces), we shall assume now that there exists $C > 0$ and
$N \in \NN$ such that
\begin{equation}\label{eq.diophantine}
	| m_1 \alpha_1 + m_2 \alpha_2 + \ldots + m_n \alpha_n |^{-1}
	\le C(|m_1| + |m_2| + \ldots |m_n|)^{N},
\end{equation}
for any $m_1, \ldots, m_n \in \ZZ$, not all zero.  (When $n = 2$, this
can be achieved by choosing $\alpha_2/\alpha_1$ to be an irrational
algebraic number, for example.)

Trivialize the normal bundle to $F$ using the standard metric on
$M$. Let $s(\ZZ^n)$ be the space of rapidly decreasing functions on
$\ZZ^n$.  The Fourier transform then establishes isomorphisms
$$
	\Omega^{k,h}(M,F) \cong s(\ZZ^n) \otimes
	\Lambda^h\CC^{n-1}\,,\quad \text{for }\;k = 0,1,
$$
and $\Omega^{k,h}(M,F) = 0$ otherwise. Under these isomorphisms, the
differential
$$
	d_F : \Omega^{0,h}(M,F) \to \Omega^{1,h}(M,F)
$$ 
becomes multiplication by $m_1 \alpha_1 + m_2 \alpha_2 + \ldots + m_n
\alpha_n$. The assumption of Equation \eqref{eq.diophantine} then
implies $\cohom^{k,h}(M,F) \cong \Lambda^h\CC^{n-1}$, for $k = 0, 1$,
and $\cohom^{k,h}(M,F) = 0$ otherwise. Thus $\cohom^{k,h}(M,F) \cong
\Lambda^k \CC \otimes \Lambda^h \CC^{n-1}$, for any $k$ and $h$.

Putting all thes calculations together we obtain
\begin{multline}
	\cohom^{k,h}(\S^*F \times \S^1, F_1) \cong (\Lambda^k \CC
	\oplus \Lambda^{k-1}\CC) \otimes \Lambda^h \CC^{n-1} \otimes
	\CC^{\{\pm \}} \\ \cong \Lambda^k \CC^2 \otimes \Lambda^h
	\CC^{n-1} \otimes \CC^{\{\pm \}}.
\end{multline}
\end{example}

\section{Complete symbols on foliations\label{Sec.CSF}}

We shall use the results of the previous sections to study the
Hochschild homology of the algebra $\alge(M,F)$ of complete symbols of
longitudinal, classical pseudodifferential operators along the leaves
of a foliation $(M,F)$. We begin by defining the algebra $\alge(M,F)$.
We assume $M$ to be compact for simplicity. Most of the following
results and constructions work for $M$ non-compact, but become
notationally more complicated. In particular, the main computations of
Hochschild homology, Theorems \ref{theorem.comp.hom} and
\ref{theorem.traces} remain true by considering compactly supported
cohomology groups.
%
In this section $n$ does not denote the dimension of $M$. 

%
If $(M,F)$ is the foliation defined by the fibers of a fibration $M
\to B$, then $\psi^\infty(M,F)$ denotes the space of smooth families
of pseudodifferential operators along the fibers of $M \to B$ and we
define $\CAA : = \psi^\infty(M,F) / \psi^{-\infty}(M,F)$.

To construct the algebra $\alge(M,F)$ in general, consider a covering
$M = \cup U_\alpha$ of $M$ with distinguished open subsets. Then
\begin{equation}\label{eq.def.AMF}
	\alge(M,F) = \sum \alge(U_\alpha,F\vert_{U_\alpha}),
\end{equation}
where the sum is taken in the space $\prod_L
\psi^\infty(L)/\psi^{-\infty}(L)$, with $L$ ranging through all leaves
of $(M,F)$. Note that it still makes sense to talk about complete
symbols of order (at most) $m$ in $\alge(M,F)$, which provides us with
a natural filtration $F_m\alge(M,F)$ of $\alge(M,F)$.
(In fact, we can define the algebra $\psi^\infty(M,F)$ similarly
\cite{ConnesIntegration, MooreSchochet}, and then $F_m \alge(M,F) =
\psi^m(M,F)/\psi^{-\infty}(M,F)$.)

Let $S^m(F^*)$ be the space of classical (compactly supported in the
base variable) symbols on the vector bundle $F^*$.  Using standard
procedures, one can define a quantization map
\cite{NistorWeinsteinXu}:
\begin{equation}
	q: S^\infty(F^*) :=\cup_{m\in \Z}S^m(F^*) \longrightarrow
	\psi^\infty(M,F),
\end{equation} 
which maps the subspace $S^m(F^*)$ of classical symbols of order $m$
to $F_m \alge(M,F)$ and satisfies $\sigma_m(q(a)) \in a +
S^{m-1}(F^*)$ if $a \in S^m(F^*)$. (In fact one could define a
quantization map $q: S^\infty(F^*) :=\cup_{m\in \Z}S^m(F^*) \to
\psi^\infty(M,F)$ descending to our quantization map, but we shall not
need this.) We can construct $q$ using a covering of $M$ by
distinguished open sets and a partition of unity. Or one can use the
results of \cite{NistorWeinsteinXu}. Our quantization map induces a
filtration preserving bijection
\begin{equation}
	S^{\infty}(F^*)/S^{-\infty}(F^*) \longrightarrow \alge(M,F)
\end{equation}

Denote now by $\GR$ the holonomy Lie groupoid associated with the
foliation $(M,F)$. The algebra $\alge(M,F)$ then coincides with the
algebra of complete symbols on $\GR$ as defined in
\cite{NistorWeinsteinXu}. Note that $\alge(M,F)$ is not a topological
algebra but it satisfies the axioms of a topologically filtered
algebra, see \cite[Proposition 3]{BenameurNistor1}. Recall that an
algebra $\alge$ with a given topology, is a topologically filtered
algebra if there exists an increasing multi-filtration $F_{n,l}^m
\alge \subset \alge$,
\begin{equation*}
 	F_{n,l}^m \alge \subset F_{n',l'}^{m'}\alge, \quad \text{if }
 	n \le n',\, l\le l',\, \text{ and } m \le m',
\end{equation*}
by closed, complemented subspaces, satisfying the following
properties:
\begin{enumerate}
\item\ $\alge = \displaystyle{ \cup_{n,l,m} } F^m_{n,l} \alge$;
\item\ The union $\alge_n := \displaystyle{ \cup_{m,l} F_{n,l}^m
\alge}$ is a closed subspace such that
\begin{equation*}
 	F_{n,l}^m \alge = \alge_n \cap \big ( \cup_{j} F_{j,l}^m \alge
 	\big );
\end{equation*}
\item\ Multiplication maps $F_{n,l}^m \alge \otimes F_{n',l'}^m \alge$
to $F_{n + n',l+l'}^m \alge$;
\item\ The maps
$$
 	F_{n,l}^m\alge /F_{n-j,l}^m \alge \otimes F_{n',l'}^m \alge
 	/F_{n'-j,l'}^m \alge \longrightarrow F_{n+n',l+l'}^m\alge
 	/F_{n + n' -j, l+l'}^m \alge
$$
induced by multiplication are continuous;\item\ The quotient
$F_{n,l}^m \alge/ F_{n-j,l}^m\alge$ is a nuclear Frechet space in the
induced topology;
\item\ The natural map
$$
 	F_{n,l}^m \alge \longrightarrow
 	\displaystyle{\lim_{\leftarrow}}\, F_{n,l}^m / F_{n-j,l}^m
 	\alge, \quad j \to \infty
$$
is a homeomorphism; and
\item\ The topology on $\alge$ is the strict inductive limit of the
subspaces $F_{n,n}^n\alge$, as $n \to \infty$ (recall that
$F_{n,n}^n\alge$ is assumed to be closed in $F_{n+1,n+1}^{n+1}\alge$).
\end{enumerate}

(The above definition corrects a typo in \cite{BenameurNistor1}, where
$n-1$ was written instead of $n-j$ in condition (4) of the above
definition.)

We have

\begin{proposition}\ Let $(M,F)$ be a smooth, compact foliated manifold.
The algebra $\CAA$ is a topologically filtered algebra such that
$F^m_{n,l}\CAA = F_n\CAA := \psi^n(M,F)/\psi^{-\infty}(M,F),$ is, in
particular, independent of $l$ and $m$.
\end{proposition}

\begin{proof}\
The algebra $\alge(M,F)$ coincides with the algebra of complete
symbols $\alge(\GR)$ on the holonomy Lie groupoid $\GR$ as studied in
\cite{BenameurNistor1} and \cite{BenameurNistor2}. Thus we get the
result by applying Proposition 3 in \cite{BenameurNistor1}.
\end{proof}

The Hochschild, cyclic, and periodic cyclic homology of the algebra
$\CAA$ must be defined by taking into account the fact that it is a
topologically filtered algebra. This is done in \cite{BenameurNistor1}
and also in \cite{BenameurNistor2}.  Fix a metric on $F$ and let $P$
be a pseudodifferential operator of order one such that $\sigma_1(P)
\equiv r$ (modulo lower order symbols), where $r \in \CI(F^*)$ is a
distance function to the origin as defined in the previous
section. The graded algebra $Gr(\CAA)$ associated to $\CAA$ is
commutative, more precisely
$$
	Gr(\CAA) \simeq \CI(\S^*F) \otimes \CC[r,r^{-1}],
$$
with grading given by the powers of $r$.

The tensor products appearing in the Hochschild complex are completed
projective tensor products such that $F_k\mathcal
H_n(\Alg)/F_{k+1}\mathcal H_n(\Alg)$ is a direct sum of spaces
isomorphic to $C^\infty( \S^*F \times \S^*F \times \ldots \times
\S^*F)$ and such that the natural map
$$
	F_k\mathcal H_n(\Alg)\longrightarrow \lim_{\leftarrow}
	F_k\mathcal H_n(\Alg)/F_{k-j}\mathcal H_n(\Alg),\quad j \to
	\infty,
$$
is an isomorphism. (This last property together with $\mathcal
H_n(\Alg) = \cup_k F_k\mathcal H_n(\Alg)$ give, by definition, the
asymptotic completeness of our Hochschild complex, see
\cite{BenameurNistor1}).

Periodic cyclic homology for algebras of complete symbols associated
with almost differentiable groupoids was computed in
\cite{BenameurNistor1}. These results include the case of the holonomy
groupoid considered in the present paper. For the sake of completeness,
let us state the explicit result for foliations.

\begin{theorem}\label{theorem.prev}\
Let $(M,F)$ be a smooth, compact foliated manifold as before. Then the
periodic cyclic homology of the algebra $\CAA$ of complete
longitudinal symbols on $(M,F)$ is given by:
$$
	\Hp_k(\CAA) \simeq \bigoplus_{j\in \Z} \cohom^{k+2j}(\S^*F
	\times \S^1), \quad k=0,1.
$$
In the same way, the periodic cyclic homology of the algebra
$\alge_0(M,F)$ of longitudinal complete symbols of order $\leq 0$ is:
$$
	\Hp_k(\alge_0(M,F)) \simeq \bigoplus_{j\in \Z}
	\cohom^{k+2j}(\S^*F), \quad k=0, 1.
$$
\end{theorem}

\section{Homology of complete symbols\label{Sec.PSDO}}

We now return to the study of the Hochschild homology of $\alge(M,F)$.
Recall that $(M,F)$ is a foliated smooth compact manifold with
$\dim(M)=n$ and $\dim(F)=p$. The codimension of the foliation will be
denoted by $q$ so $n=p+q$.

The canonical filtration of the Hochschild complex defined above
(following \cite{BenameurNistor1}) gives rise to a spectral sequence
$\EH_{k,h}^r$, by general results about filtered complexes.  This
spectral sequence has the $\EH^1$-term given by
$$
	\EH^1_{k,h}=\HH_{k+h}(Gr(\CAA))_k,
$$
by \cite[Lemma 1]{BenameurNistor1}.  The Hochschild homology of
$Gr(\CAA)$ is identified using a combination of the
Hochschild-Kostand-Rosenberg (HKR) isomorphism and a result of Connes,
which is the analog of the HKR-isomorphism for algebras of smooth
functions.  We denote by $\Omega^{r}(F^*\smallsetminus M)_s$ the set
of differential $r-$forms on the manifold $F^*\smallsetminus M$ which
are positively $s$-homogeneous in the radial direction. Then we have:
$$
	\Hd_{l}(Gr(\Alg))_{d} \simeq \Omega^{l}(\S^*F)r^d \oplus
	\Omega^{l-1}(\S^*F)r^{d-1}dr,
$$
the isomorphism being obtained via the
Hochschild-Kostant-Rosenberg-Connes map
$$
	\chi(a_0, \ldots , a_l)={\frac{1}{l!}} a_0 da_1 \ldots da_l.
$$

Let $X := F^* \smallsetminus M$, as above. It will be convenient to
identify
$$
	\Omega^{l}(\S^*F)r^d \oplus \Omega^{l-1}(\S^*F)r^{d-1}dr
$$
with the subspace $\Omega^{l}(X)_d \subset \Omega^{l}(X)$ consisting
of $d$-homogeneous $l$-forms on the manifold $X = F^* \smallsetminus
M$.  Also, we endow $X$ with the foliation $\maF \subset TX$ whose
leaves are the cotangent bundles to the leaves of $(M,F)$ with the
zero section removed, as we did in Sections \ref{Sec.deRham} and
\ref{Sec.Can.Hom}. More precisely, if $\pi: F^* \to M$ is the
projection and $T_vF^* = \ker \pi_*$ is the vertical tangent bundle to
the fibration $F^* \to M$, then $\maF$ is the restriction to $X$ of
the bundle $\pi^*(F) + T_v F^*$.

Recall that $X = F^* \smallsetminus M$ admits a Poisson structure
induced by the natural symplectic structure of the leaves of
$\maF$. Moreover $(X, \maF)$ is then a conic symplectic foliation in
the sense of Definition \ref{conic}.  We introduced in Section
\ref{Sec.Can.Hom}, a Poisson differential
$\delta=\delta_{\maF}+\delta_{-2,1} : \Omega^{l}(X) \rightarrow
\Omega^{l-1}(X)$, such that $\delta \left(\Omega^{l}(X)_k \right)
\subset \Omega^{l-1}(X)_{k-1}$. We denote as in Section
\ref{Sec.Can.Hom} by
$$
	\cohom_{l}^\delta(X)_d = {\frac{\ker (\delta:\Omega^l(X)_d \to
	\Omega^{l-1}(X)_{d-1})}{ \delta (\Omega^{l+1}(X)_{d+1}) }}
$$ 
the homogeneous Poisson homology groups of $X=F^*\smallsetminus M$.

\begin{proposition}\label{prop.E1}\
Let $\chi: \HH_{l}(Gr(A))_d \rightarrow \Omega^{l}(F^*\smallsetminus
M)_d$ be the HKR isomorphism, and let $d_1 : E^1_{k,h}\rightarrow
E^1_{k-1,h}$ be the first differential of the spectral sequence
associated to $\CAA$ as in \cite{BenameurNistor1}.  Then
$$
	\chi \circ d_1 \circ \chi^{-1}=-\sqrt{-1} \delta,
$$
and hence $\EH_{k,h}^2 \simeq \cohom^\delta_{k+h}(F^*\smallsetminus
M)_k$.
\end{proposition}

\begin{proof}\ We apply Theorem 3.1.1 in \cite[page 107]{Brylinski}. 
More precisely, if $\sigma\in {\mathcal{A}}_m /{\mathcal{A}}_{m-1}$,
$\sigma'\in {\mathcal{A}}_q/{\mathcal{A}}_{q-1}$ then there exist
$$
	P \in \psi^m(M,F) \; \text{ and } \; P'\in \psi^q(M,F)
$$ 
such that $[P]=\sigma$ and $[P']=\sigma'$ (with obvious
notations). One needs the expansion of $P\circ P'-P'\circ P$ into
homogeneous terms. Since the quatization map can be chosen with values
operators with support small enough, the symbol expansion of the
commutator is obtained in the same way as in the classical case, see
\cite{Gilkey} and \cite{ConnesIntegration, ConnesFOL} for the
corresponding results for foliations.
\end{proof}

Let $F_1 =\pi^*(F) + T_v(\S^*F) \subset T(\S^*F)$ be the integrable
sub-bundle defined using $F$, as above, but for the cosphere
bundle. (By abuse of notation, we shall sometimes denote $F_1$ also
the integrable sub-bundles defined similarly by $F$ on the fibrations
$F^* \to M$, on $F^* \smallsetminus M \to M$, on $\S^*F \to M$, or on
$\S^*F \times \S^1 \to M$. So on $F^*\smallsetminus M$ for instance
$F_1$ coincides with the foliation $\maF$ defined in the previous
sections.)  Now we gather the results of the previous sections and
deduce the following theorem.

\begin{theorem}\label{theorem.comp.hom}\
Let $(M,F)$ be a smooth, compact foliated manifold and denote by $p$
the dimension of the leaves of $(M,F)$. Let $(\EH^r,d^r)_{r\geq 1}$ be
the spectral sequence associated to the canonical filtration of the
Hochschild complex of the algebra $\Alg (M,F) := \psi^\infty(M,F)/
\psi^{-\infty}(M,F)$, then this spectral sequence converges to the
Hochschild homology of $\CAA$ and its $E^2$ term is given by
$$
	\EH^2_{k,h} \simeq\cohom^{p - k, h - p} (\S^*F \times \S^1,
	F_1).
$$
\end{theorem}

\begin{proof}\
First we have an isomorphism $\EH^2_{k,h} \cong
\cohom^\delta_{k+h}(F^*\smallsetminus M)_k$ given by Proposition
\ref{prop.E1}. By Theorem \ref{hom.can} we have
$$
	\cohom^\delta_{k+h} (F^*\smallsetminus M, F_1)_k \cong
	\cohom^{p-k,h-p}(F^*\smallsetminus M, F_1)_0.
$$
But this last group coincides with $\cohom^{p - k, h - p} (\S^*F
\times \S^1, F_1)$, as we have already checked in the proof of
Corollary \ref{cor.hom.can}.

The convergence of the spectral sequence is then a consequence of
\cite[Lemma 3]{BenameurNistor1} by taking $a=2$ in that lemma.
\end{proof}

\begin{theorem}\label{theorem.traces}\
Let $(M,F)$ be a smooth, compact foliated manifold with $\dim(F)=p$.
Then the space $\HH_0(\CAA)$ is given by:
$$
	\HH_0(\CAA) \simeq \cohom^{2p,0}(\S^*F \times \S^1,F_1).
$$
Moreover, when $p\geq 2$, we have
$$
	\HH_0(\CAA) \simeq \cohom^{p,0}(M,F).
$$
\end{theorem}

\begin{proof}\ 
Using the previous results, we need to show that the differentials
$d^r$ coming into and out of $\EH^r_{k,h}$ are trivial if $k+h=0$.
But the $E^2_{k,h}$ term vanishes unless $-p \le k \le p$ and $p \le h
\le p + q$ where $q$ is the codimension of the foliation. Thus the
only term $E^2_{-k,k}$ that may be different from $0$ is
$E^2_{-p,p}$. All differentials coming into and out of $E^2_{-p,p}$
are seen to vanish because of the geometry of this spectral
sequence. More precisely, recall on that $d^2:\EH^2_{k,h} \to
\EH^2_{k-2, h+1}$ while $\EH^2_{k,h}$ is only non trivial when $-p\leq
k\leq p$ and $p\leq h \leq p+q$. Thus
$$
	d^2|_{\EH^2_{-p,p}} \text{ is trivial},
$$
and the range of $d^2$ does not intersect $\EH^2_{-p,p}$.  In the same
way,
$$
	d^r:\EH^r_{k,h} \to \EH^r_{k-r,h+r-1},
$$ 
thus it is trivial when $k=-p$ and $h=p$, and its range is never of
the form $\EH^r_{-p,p}$.  A recursive argument then finishes the
proof. By Corollary \ref{Cor.2}(ii), when $2p-1\geq p+1$, we have
$$
	\cohom^{2p,0}(\S^*F\times \S^1, F_1) \simeq \cohom^{p,0}(M,F).
$$
Thus the proof is complete.
\end{proof}

\begin{remark}\
The above restriction $p\geq 2$ corresponds to the connectedness of
the total manifold of the bundle $\S^*F$ and is similar to the
restriction on the uniqueness of the Wodzicki trace in the non
foliated situation.
\end{remark}

\begin{remark}\
A similar argument enables to obtain that
$$
	\HH_{2p+q} (\maA(M,F)) \simeq \cohom^{0,q}(M,F).
$$
This follows from the Gysin exact sequence, Theorem \ref{Gysin}.
\end{remark}

It was proved in \cite{BenameurNistor2} that the above spectral
sequence collapses at $E^2$ when the given foliation is a smooth
fibration. It would be interesting to establish this result in
general, because of the following corollary.

\begin{corollary}\label{cor.comp.HH}\
Let $(M,F)$ be a smooth, compact foliated manifold with
$\dim(F)=p$. Assume that the spectral sequence associated with
Hochschild homology collapses at $E^2$, then the Hochschild homology
of the algebra of complete longitudinal pseudodifferential symbols on
$(M,F)$ is given by:
$$
	\HH_k(\CAA) \simeq \bigoplus_{j=0}^q \cohom^{2p + j - k,
	j}(\S^*F\times \S^1, F_1).
$$
\end{corollary}

\begin{proof}\ We have 
\begin{multline*}
	\HH_k(\CAA)\simeq \oplus_{l} \EH^2_{k - l,l} \simeq \oplus_l
	\cohom^{p - k + l , l - p} (F^*\smallsetminus M, F_1)_0 \\
	\simeq \oplus_{j=0}^q \cohom^{2p + j - k, j}
	(F^*\smallsetminus M, F_1)_0.
\end{multline*}
Using Corollary \ref{cor.hom.can}, it remains to show the convergence
of the above spectral sequence, but this was checked in Theorem
\ref{theorem.comp.hom}.
\end{proof}

Let us now formulate the corresponding results for Hochschild
cohomology. First, to define Hochschild cohomology, we just dualize
the constructions (inductive and projective limits, but keeping the
projective limits first) used to define the Hochschild homology
complex in \cite{BenameurNistor1}.  In particular, {\em all cocycles
$\phi$ in the definition of Hochschild complex are such that
$\phi(a_0,\ldots,a_k) = 0$ if the sum of the orders of $a_0, a_1,
\ldots , a_k$ is less than some fixed number $N$, that is fixed for
each $\phi$.} The same theorems on the convergence of the associated
spectral sequences then hold for Hochschild cohomology (with the same
proof).

\begin{theorem}\ 
Let $(M,F)$ be a smooth, compact foliated manifold with
$\dim(F)=p$. Let $F_1$ be the foliation of $\S^*F\times \S^1$ induced
by $F$ and of the same codimension as $F$, as above.  Then the
spectral sequence $\EH_r^{k,h}$ associated with the Hochschild
cohomology of $\CAA$ converges to Hochschild cohomology and has the
$\EH_2$-term given by $ \EH_2^{k,h} \simeq \cohom_{p-k,h-p}(\S^*F
\times \S^1, F_1).$

In particular, $\HH^0(\CAA) \simeq \cohom_{2p,0}(\S^*F \times \S^1,
F_1)$.
\end{theorem}

Thus traces are constructed out of $(2p,0)$-invariant currents on
$(\S^*F\times \S^1,F_1)$.  But for $p\geq 2$, we have a homological
identification, similar to the cohomological one obtained in Theorem
\ref{Gysin}:
$$
	\cohom_{2p,0}(\S^*F \times \S^1, F_1) \simeq \cohom_{p,0}(M,
	F).
$$
The space $\cohom_{p,0}(M,F)$ is the space of closed holonomy
invariant $p$-currents, see \cite{AbouqatebKacimi}.

\begin{example} 
Let us take a closer look at the foliation of the Example
\ref{ex.torus}. By duality, we obtain
\begin{equation}
	\cohom_{k,h}(\S^*F \times \S^1, F_1) \cong \Lambda^k\CC^2
	\otimes \Lambda^h\CC^{n-1} \otimes \CC^{\{\pm \}}.
\end{equation}
This gives us a canonical basis, $\tau_{+}$ and $\tau_{-}$, for
$\HH^0(\CAA)$.  It also gives that the dimension of $\HH^l(\CAA)$ is
at most the dimension of the space $\Lambda^l(\CC^{n+1})\otimes
\CC^2$, and that these dimensions are equal if, and only if, the
spectral sequence $\EH_r$ collapses at $\EH_2$.

Chose a subtorus of codimension $1$ in $M$ that is transverse to the
foliation. This gives rise to $n-1$ one-parameter groups of
automorphisms of $\CAA$, and hence to $n-1$-derivations $\delta_1,
\ldots, \delta_{n-1}$ of this algebra.  Let $\delta$ be the derivation
given by translation along the leaves of the foliation and
$\delta_r(T) = [\log Q,T]$, where $Q$ is a positive operator of order
$1$ with principal symbol $r$, as in \cite{MelroseNistor}. Then each
of these derivations acts on the Hochschild complex of $\CAA$.

If $D$ is a derivation and $\phi$ is a $l$-cocycle on $\CAA$, then
\begin{equation}\label{eq.cocycles}
	(i_D\phi)(a_0,a_1, \ldots , a_{l+1}) := \phi(a_0D(a_1), a_2,
	\ldots, a_{l+1})
\end{equation}
will be a $(l+1)$-cocycle on $\CAA$.  We have $i_D i_{D'} =
-i_{D'}i_{D}$, for all $D,D' \in \{\delta_1, \ldots, \delta_{n-1},
\delta, \delta_r\}$, because all these derivations commute.

A counting argument gives then that there are as many $l$-cocycles of
the form $i_{D_1}\ldots i_{D_l}\tau_{\pm }$, with $D_1, \ldots, D_l$
distinct elements in the set $\{\delta_1, \ldots, \delta_{n-1},
\delta, \delta_r\}$ as the maximum possible dimension of $\HH^l(\CAA)$
established above (that is, the dimension of $\Lambda^l\CC^{n+1}
\otimes \CC^2$.

The algebra $\CAA$ splits canonically as a direct sum
$$
	\CAA \cong \CAA_+ \oplus \CAA_-,
$$
because the contangent sphere bundle $\S^*F$ is disconnected.  We next
use the inclusion $\CI(M) \subset \CAA$ and the fact that $i_{D_r}$
induces a morphism of the Hochschild complexes to prove that all the
cocycles in Equation \eqref{eq.cocycles} are distinct. This shows that
the spectral sequence $\EH_r$ degenerates at $\EH_2$.  It also gives
an explicit determination of a basis of the groups $\HH^l(\CAA)$ for
this foliation. In particular
$$
	\HH^l(\CAA) \cong \Lambda^l \CC^{n+1} \otimes \CC^{\{\pm\}}.
$$
\end{example}

\end{document}